\documentclass[12pt,a4paper]{article}

\usepackage[latin1]{inputenc}
\usepackage{amsmath}
\usepackage{amsfonts}
\usepackage{amssymb}

\author{Pedro Caro \footnote{Department of Mathematics, Universidad Aut\'onoma de Madrid, 28049 Madrid, Spain}}
\title{Stable Determination of the Electromagnetic Coefficients by Boundary Measurements}
\date{March 17, 2010}

\newcommand{\inte}{\Omega}
\newcommand{\curl}{\nabla \! \times \!}
\newcommand{\diver}{\nabla \cdot}
\newcommand{\Durl}{D \times \!}
\newcommand{\Diver}{D \cdot}
\newcommand{\bou}{\partial\Omega} 
\newcommand{\Div}{\textrm{Div}} 

\newcommand{\frv}[4]{\left(\begin{array}{c} 
#1\\ #2\\
\hline
#3\\ #4\\
\end{array}\right)}
\newcommand{\fcv}[4]{\left(\begin{array}{c c|c c} 
#1 & #2 & #3 & #4
\end{array}\right)}
\newcommand{\fbM}[4]{\left(\begin{array}{c|c} 
#1 & #2\\
\hline
#3 & #4
\end{array}\right)}
\newcommand{\subM}[4]{\begin{array}{c c} 
#1 & #2\\
#3 & #4
\end{array}}
\newcommand{\inner}[2]{#1 \! \cdot \! #2} 
\newcommand{\Inner}[2]{\left( #1 \middle| #2 \right)} 
\newcommand{\norm}[3]{\left\|#1\right\|^{#2}_{#3}}  
\newcommand{\cross}[2]{#1 \! \times \! #2} 
\newcommand{\dual}[2]{\left\langle #1 \middle| #2 \right\rangle} 
\newcommand{\Hdiv}[1]{H^{#1}(\inte; \textrm{div})} 
\newcommand{\Hcurl}[1]{H^{#1}(\inte; \textrm{curl})} 

\newtheorem{theorem}{Theorem}
\newtheorem{definition}{Definition}
\newtheorem{lemma}[definition]{Lemma}
\newtheorem{proposition}[definition]{Proposition}
\newtheorem{corollary}{Corollary}
\newenvironment{proof}{\textbf{Proof:}} {\hspace{\stretch{1}}$\Box$}

\begin{document}

\maketitle

\begin{abstract}
The goal of this paper is to prove a stable determination of the coefficients for the time-harmonic Maxwell equations, in a Lipschitz domain, by boundary measurements.
\end{abstract}

\section*{Introduction}
Let $ \inte $ be a bounded Lipschitz domain in the three-dimensional euclidean space. Assume the medium, modeled by $ \inte $, to be non-homogeneous and isotropic. Suppose the electromagnetic properties of $ \inte $ to be described by the electric permittivity $ \varepsilon $, the magnetic permeability $ \mu $ and the electric conductivity $ \sigma $. Let $ E, H $ denote the electric and magnetic fields, respectively. The time-harmonic Maxwell equations at frequency $ \omega $ read
\[ \curl H + i \omega \varepsilon E = \sigma E, \qquad \curl E - i \omega \mu H = 0, \]
whenever the total electric current density is given by $ \sigma E $.
Writing $ \gamma = \varepsilon + i \sigma/ \omega $, the time-harmonic Maxwell equations can be expressed as
\begin{equation}\label{ME_without_sources}
\left\{ 
\begin{array}{l}
\curl H + i \omega \gamma E = 0\\
\curl E - i \omega \mu H = 0.
\end{array}
\right.
\end{equation}
It is known that this system, complemented with a suitable prescribed data on the boundary, is well-posed for some $ \omega $'s. In fact, we have the following result.

\begin{theorem}\label{th:direct_problem}\rm
Let $ \inte $ be a bounded Lipschitz domain and $ \mu, \varepsilon, \sigma \in L^\infty(\inte) $ satisfying 
\[ \mu \geq \mu' > 0, \quad \varepsilon \geq \varepsilon' > 0, \quad \sigma \geq 0; \]
a. e. in $ \inte $, with $ \mu', \varepsilon' $ positive constants. Given $ T \in TH(\bou) $, the problem of finding $ E, H \in \Hcurl{} $ solving (\ref{ME_without_sources}) in $ \inte $ and satisfying either $ \cross{N}{E} = T $ or $ \cross{N}{H} = T $ is well-posed for any $ \omega \in \mathbb{C} \setminus \{ 0 \} $ except for a subset of
\[ \{ \omega \in \mathbb{C} : -\norm{\sigma / \varepsilon}{}{L^\infty(\inte)} \leq \mathrm{Im}\, \omega \leq 0 \} \]
with no accumulation point in $ \mathbb{C} \setminus \{0\} $.
\end{theorem}
Here $ N $ stands for the unit vector field normal to $ \bou $, the boundary of $ \inte $. Precise definitions of the spaces $ \Hcurl{} $ and $ TH(\bou) $ are given in Subsection \ref{sec:non-standar_SovBes}. The frequencies $ \omega $ for which the direct problem is not well-posed are called \emph{resonant frequencies}.

A proof of this result in the case where $ \bou $ is $ C^2 $ can be found in \cite{SIsCh}. The well-posedness of the same problem in the context of non-smooth domains is studied in \cite{Le} but the result stated there assumes $ \sigma $ to be zero. Finally, a proof of the precise statement of Theorem \ref{th:direct_problem} can be given using the compactness result stated in \cite{W} and a lemma due to Peetre (see \cite{Pe} and \cite{LMa}).

Theorem \ref{th:direct_problem} allows us to model the boundary measurements by means of the admittance or the impedance maps. The \emph{admittance map} is defined as
\[ \Lambda^{ad} : T \in TH(\bou) \longmapsto \cross{N}{H} \in TH(\bou), \]
where $ E, H $ is the solution for (\ref{ME_without_sources}) with $ \cross{N}{E} = T $. The \emph{impedance map} is defined as
\[ \Lambda^{im} : T \in TH(\bou) \longmapsto \cross{N}{E} \in TH(\bou), \]
where $ E, H $ is the solution for (\ref{ME_without_sources}) with $ \cross{N}{H} = T $.
It is a consequence of Theorem \ref{th:direct_problem} that these maps can only be used out of resonant frequencies.

With one of these maps at hand, one can set the problem of determining the coefficients of equations in (\ref{ME_without_sources}). This inverse problem was initially proposed by Somersalo \textit{et al} in \cite{SIsCh}. In these notes we prove not only determination of the coefficients but also stable determination.

When trying to determine in a stable manner the coefficients by one of these maps, we have to face the problem of choosing $ \omega > 0 $ to be non-resonant for the class of coefficients to determine. How can we manage to solve this if the position of the resonant frequencies depends on the unknown coefficients? As it was pointed out in \cite{OPS2}, the same happens in practice. How can we know when our choice of frequency is close to a resonant frequency?

In order to avoid this problem we shall model the boundary measurements by the Cauchy data set. Given a frequency $ \omega > 0 $, the \emph{Cauchy data set} is defined as follows: $ (T, S) \in C(\mu, \gamma) $ if and only if $ (T, S) \in (TH(\bou))^2 $ and there exists a pair $ (E, H) \in (\Hcurl{})^2 $ solution of (\ref{ME_without_sources}) satisfying $ \cross{N}{E} = T $ and $ \cross{N}{H} = S $.

This way of modeling the boundary measurements has been used successfully in \cite{BuU}, \cite{SaTz} and \cite{SaTz2}. However, as far as the author knows Cauchy data sets have not yet been used to establish any result of \emph{stability}.

In order to measure the proximity of the Cauchy data sets associated to given pairs $ \mu_1, \gamma_1 $ and $ \mu_2, \gamma_2 $ we introduce a pseudo-metric distance.

\begin{definition}\label{def:Cdistance}\rm
Let  $ \mu_1, \gamma_1 $ and $ \mu_2, \gamma_2 $ be two pairs of coefficients. Consider $ \omega $ a positive frequency and let $ C_j $ denote $ C(\mu_j, \gamma_j) $. Let us define the pseudo-metric distance between the Cauchy data sets $ C_1 $ and $ C_2 $ as
\begin{equation*}
\delta_C (C_1, C_2) = \max_{j \neq k} \sup_{ \substack{(T_k, S_k) \in C_k \\ \norm{T_k}{}{TH(\bou))} = 1 }} \inf_{(T_j, S_j) \in C_j} \norm{(T_j, S_j) - (T_k, S_k)}{}{(TH(\bou))^2}.
\end{equation*}
\end{definition}
We say that $ \delta_C $ is a pseudo-metric distance because if $ \delta_C(C_1, C_2) = 0 $, we can only ensure that $ \overline{C_1} = \overline{C_2} $.

The definition of $ \delta_C $ is inspired in the Hausdorff distance. Unlike the latter distance the former one is comparable to $ \norm{\Lambda^{ad}_1 - \Lambda^{ad}_2}{}{} $ when $ \omega $ is a non-resonant frequency for $ \mu_j, \gamma_j $ with $ j =1,2 $. Here $ \norm{\centerdot}{}{} $ denotes the operator norm for linear operators on $ TH(\bou) $.

Our result requires certain stability of the problem on the boundary and since this has not been proven yet, we shall introduce some definitions.
\begin{definition}\label{def:admissible}\rm
Given two constants $ M, s $ such that $ 0 < M $, $ 0 < s < 1/2 $, we shall say that the pair of coefficients $ \mu, \gamma $ is \emph{admissible} if they satisfy the following conditions.
\begin{itemize}
\item[(i)] \emph{Uniform ellipticity condition.} The coefficients $ \gamma, \mu \in C^{1,1}(\overline{\inte}) $ satisfy
\[ M^{-1} \leq \mathrm{Re}\,\gamma(x) \qquad M^{-1} \leq \mu(x); \]
for any $ x \in \inte $.
\item[(ii)] \emph{A priori bound on the boundary.} The following a priori bound holds on the boundary
\begin{equation*}
\norm{\gamma}{}{C^{0,1}(\bou)} + \norm{\mu}{}{C^{0,1}(\bou)} < M.
\end{equation*}
\item[(iii)] \emph{A priori bound in the interior.} The following a priori bounds hold in the interior
\begin{align*}
\norm{\gamma}{}{W^{2,\infty}(\inte)} + \norm{\mu}{}{W^{2,\infty}(\inte)} \leq M,& &
\norm{\gamma}{}{H^{2 + s}(\inte)} + \norm{\mu}{}{H^{2 + s}(\inte)} \leq M.
\end{align*}
\end{itemize}
\end{definition}

\begin{definition}\rm
Let $ M, s $ be the constants given in Definition \ref{def:admissible} and let $ \omega $ be a positive frequency. We shall say that a pair $ \mu, \gamma $ is in the \emph{class of $ B $-stable coefficients on the boundary at frequency $ \omega $} if $ \mu, \gamma $ is an admissible pair and there exists a modulus of continuity $ B $ such that, for any other admissible pair $ \tilde{\mu}, \tilde{\gamma} $, one has
\begin{gather*}
\norm{\gamma - \tilde{\gamma}}{}{C^{0,1}(\bou)} + \norm{\mu - \tilde{\mu}}{}{C^{0,1}(\bou)} \leq B \left( \delta_C(C, \tilde{C}) \right),\\
\norm{\nabla (\gamma - \tilde{\gamma})}{}{L^\infty(\bou; \mathbb{C}^3)} + \norm{\nabla (\mu - \tilde{\mu})}{}{L^\infty(\bou; \mathbb{C}^3)} \leq B \left( \delta_C(C, \tilde{C}) \right).
\end{gather*}
\end{definition}
Here $ C, \tilde{C} $ are the Cauchy data sets associated to the pairs $ \mu, \gamma $ and $ \tilde{\mu}, \tilde{\gamma} $, respectively.

With these definitions at hand the stable determination of the electromagnetic coefficients can be stated as follows. 
\begin{theorem}\label{th:stability}\rm
Let $ \inte $ be a bounded Lipschitz domain and let $ \omega $ be a positive frequency. Consider $ \mu_1, \gamma_1 $ and $ \mu_2, \gamma_2 $ any two pairs in the class of $ B $-stable coefficients on the boundary at frequency $ \omega $, with $ B $ satisfying $ |r| \leq B(|r|) $ for all $ |r| < 1 $. Then, there exists a constant $ C = C(M) $ such that the following estimate holds
\[ \norm{\gamma_1 - \gamma_2}{}{H^1(\inte)} + \norm{\mu_1 - \mu_2}{}{H^1(\inte)} \leq C |\mathrm{log}\, B(\delta_C(C_1, C_2))|^{- \lambda}, \]
for some constant $ \lambda $ such that $ 0 < \lambda < 2/3 s $. Here $ C_1, C_2 $ are the Cauchy data sets associated to the pairs $ \mu_1, \gamma_1 $ and $ \mu_2, \gamma_2 $, respectively.
\end{theorem}
As in the inverse conductivity problem, it should be possible to prove that any admissible pair is in the class of H\"older-stable coefficients on the boundary for any frequency $ \omega $, that is, with $ B(|r|) = |r|^\alpha $ for $ 0 < \alpha < 1 $. Notice that in the conductivity case a logarithmic module of continuity, as the one in Theorem \ref{th:stability}, is optimal (see \cite{Man}).

In \cite{OPS} and \cite{OS} the coefficients were assumed to be constant on the boundary --which is actually quite natural from the point of view of applications--, in that case our result reads as follows.
\begin{corollary}\rm
Let $ \inte $ be a bounded Lipschitz domain and let $ \omega $ be a positive frequency. Consider $ \mu_1, \gamma_1 $ and $ \mu_2, \gamma_2 $ any two pairs of admissible coefficients --in the sense of Definition \ref{def:admissible}. Assume that
\[ \mu_1|_{\bou} = \mu_2|_{\bou}, \, \partial_{x^j} \mu_1|_{\bou} = \partial_{x^j} \mu_2|_{\bou},\quad \gamma_1|_{\bou} = \gamma_2|_{\bou},\, \partial_{x^j} \gamma_1|_{\bou} = \partial_{x^j} \gamma_2|_{\bou}, \]
with $ j = 1,2,3 $. Then, there exists a constant $ C = C(M) $ such that the following estimate holds
\[ \norm{\gamma_1 - \gamma_2}{}{H^1(\inte)} + \norm{\mu_1 - \mu_2}{}{H^1(\inte)} \leq C |\mathrm{log}\, \delta_C(C_1, C_2)|^{- \lambda}, \]
for some constant $ \lambda $ such that $ 0 < \lambda < 2/3 s $.
\end{corollary}

When the boundary is smooth, at least $ C^{1,1} $, it is possible to prove the existence of $ H^1 $-solutions, this allows us to describe traces not only in a weak sense but also in a strong sense. In consequence, the boundary measurements can be modeled in \emph{better} spaces (see \cite{OPS} and \cite{OS}). Furthermore, it is well known (see \cite{Sr}, \cite{Sr1} \cite{BiSo}, \cite{BiSo2} and \cite{CoD}) that Maxwell's equations may not admit $ H^1 $-solutions even with $ H^{1/2} $-boundary data, whenever the domain is neither convex nor has $ C^{1,1} $-boundary. When working on Lipschitz domain, this lack of smoothness makes necessary to introduce some non-standard Sobolev spaces.

The unique recovery of $ C^3 $-coefficients $ \gamma $ and $ \mu $ from boundary data was proved in \cite{OPS}, and later simplified in \cite{OS}, in the context of $ C^{1,1} $-domains. The regularity in our result of stability agrees with the required in \cite{OPS} and \cite{OS}, nevertheless from our proof one can state a uniqueness  result for $ C^{1,1} $-coefficients in the context of Lipschitz domains. Boundary determination results were given in \cite{McD} and \cite{JMcD} in the case that the boundary is smooth. The more general chiral media was studied in \cite{McD2}.  For a slightly more general approach and more background information, see also the review article \cite{OPS2}.

In \cite{Sk} the relation between the inverse scattering problem (for short ISP) and the inverse boundary value problem (for short IBVP) was studied. Further, Sarkola proved that the former problem can be reduced to the latter one. On the other hand, in \cite{H} H\"ahner gives a proof of the stability for the ISP when magnetic permeability $ \mu $ is constant and as a consequence he obtains stability for the IBVP for constant magnetic permeability when the domain is a ball. Thus, our corollary generalizes the stability result stated in \cite{H}, at least, for the IVBP.

The main motivation to study the inverse problem in the setting of Lipschitz domain comes from the technical hypothesis that had to be assumed in \cite{COSa} when studying local data result. Some of the ideas shown along these notes can be adapted to prove a stable determination of the coefficients from local data. This is accomplished in \cite{C3}.

The structure of the paper is as follows. There are two sections and an appendix, in the first section some preliminary details are given, while in the second one we prove the \textit{stability} of the IBVP. The appendix contains basic definitions and properties of the spaces used along these notes. In a deeper extent, the first subsection of Section \ref{sec:prelim} is dedicated to the functional spaces adapted to Maxwell's equation and their traces. The definitions of these spaces and traces have not been placed in the appendix in order for the reader to get used to them. In the second subsection of Section \ref{sec:prelim}, we transform the system (\ref{ME_without_sources}) into a Schr\"odinger-type equation. In Section \ref{sec:log-type_estimate} we perform some of the standard steps when studying this IBVP: prove a suitable estimate relating the boundary measurements with the coefficients in the interior, construct special solutions for the Maxwell's equations through the Schr\"odinger-type equation and plug these solutions into the suitable estimate. Finally, we obtain the estimate stated in Theorem \ref{th:stability} by the use of a Carleman estimate.

\paragraph{Acknowledgement.} This paper is part of the author's doctoral dissertation and it has been written under the supervision of Alberto Ruiz. The author would like to thank him for his support --\textit{gracias por dejarme elegir y corregir mis errores}. Additionally, the author would like to thank the comments of the referees. The presentation of the result has improved a lot with their suggestions. The author was economically supported by Ministerio de Ciencia e Innovaci\'on de Espa\~na, MTM2008-02568-C02-01.

\section{Preliminaries}\label{sec:prelim}
Along this section we introduce some spaces that turn to be useful in the context of Maxwell's equations and we transform (\ref{ME_without_sources}) into a Schr\"odinger-type equation.

Let us denote by $ \mathbb{E} $ the \emph{three-dimensional euclidean space} and by $ \mathcal{T} \mathbb{E} $ the \emph{module of smooth vectors fields} over the real smooth functions. The elements $ u \in \mathcal{T} \mathbb{E} $ will be called \emph{real vector fields} and can be expressed as
\begin{equation*}
u = \left(\begin{array}{c c c}
u^{(1)} & u^{(2)} & u^{(3)}
\end{array}\right)^t.
\end{equation*}
The euclidean metric induces on $ \mathbb{E} $ a \emph{volume element} $ dV $ and the \emph{euclidean distance} $ d_e $. Additionally, for any two vector fields $ u, v $, we shall denote by $ u \cdot v $ and $ u \times v $ the \emph{point-wise inner product} and the \emph{point-wise cross product} of vector fields, respectively.
On the other hand, denote
\[ \mathcal{X} \mathbb{E} = \{ u + iv : u,v \in \mathcal{T} \mathbb{E} \} .\]
The elements of $ \mathcal{X} \mathbb{E} $ will be called \emph{complex vector fields}.

Finally, the \emph{gradient}, \emph{divergence} and \emph{curl} operators will be denoted by $ \nabla $, $ \diver $ and $ \curl $, respectively.

\subsection{Non-standard Sobolev and Besov spaces}\label{sec:non-standar_SovBes}
Most of the facts collected here can be found in \cite{M}.

\begin{definition}\rm
Define the spaces
\begin{gather*}
\Hdiv{} = \{ u \in L^2(\inte; \mathbb{C}^3) : \diver u \in L^2(\inte) \},\\
\Hcurl{} = \{ v \in L^2(\inte; \mathbb{C}^3) : \curl v \in L^2(\inte; \mathbb{C}^3) \}
\end{gather*}
equipped with the graph norms
\begin{gather*}
\norm{u}{}{\Hdiv{}} = \norm{u}{}{L^2(\inte; \mathbb{C}^3)} + \norm{\diver u}{}{L^2(\inte)},\\
\norm{v}{}{\Hcurl{}} = \norm{v}{}{L^2(\inte; \mathbb{C}^3)} + \norm{\curl v}{}{L^2(\inte; \mathbb{C}^3)}.
\end{gather*}
\end{definition}

Now we define the traces of elements belonging to these spaces. For any $ u \in \Hdiv{} $ the \emph{normal trace} of $ u $, that is, $ N \cdot u $, can be defined as an element of $ B^{-1/2}(\bou) $. Namely, for any $ g \in B^{1/2}(\bou) $
\begin{equation}\label{def:normal-component}
\dual{\inner{N}{u}}{g} = \int_\inte (\diver u) \overline{f} \, dV + \int_\inte \inner{u}{\overline{\nabla f}} \, dV,
\end{equation}
where $ f \in H^1(\inte) $ and $ f|_{\bou} = g $. It is well-known that
\[ N \cdot \centerdot : \Hdiv{} \longrightarrow B^{-1/2}(\bou) \]
is a bounded and onto operator. Additionally, note that if $ u \in \Hdiv{} $, $ f \in C^{0,1}(\bou) $ and $ \tilde{f} $ is any extension of $ f $ such that $ \tilde{f} \in C^{0,1}(\overline{\inte}) $ (see \cite{St} for extensions from closed set), then
\begin{equation}\label{for:div-product}
\dual{f\inner{N}{u}}{g} = \dual{\inner{N}{u}}{g\overline{f}} = \dual{\inner{N}{(\tilde{f} u)}}{g},
\end{equation}
for any $ g \in B^{1/2}(\bou) $.

Similarly, for any $ u \in \Hcurl{} $ the \emph{tangential trace} of $ u $, that is, $ \cross{N}{u} $, can be defined as an element of $ B^{-1/2}(\bou; \mathbb{C}^3) $. Namely, for any $ w \in B^{1/2}(\bou; \mathbb{C}^3) $
\[ \dual{\cross{N}{u}}{w} = \int_\inte \inner{(\curl u)}{\overline{v}} \, dV - \int_\inte \inner{u}{(\overline{\curl v})} \, dV, \]
where $ v \in H^1(\inte; \mathbb{C}^3) $ and $ v|_{\bou} = w $. In this context, the operator
\[ \cross{N}{\centerdot} : \Hcurl{} \longrightarrow B^{-1/2}(\bou; \mathbb{C}^3) \]
is bounded but not onto. Let its range be denoted by
\[ TH(\bou) = \{ w \in B^{-1/2}(\bou; \mathbb{C}^3) : \exists u \in \Hcurl{}, \, \cross{N}{u} = w \}. \]
This vector space can be equipped with the norm
\[ \norm{w}{}{TH(\bou)} = \inf \{ \norm{u}{}{\Hcurl{}} : u \in \Hcurl{}, \, \cross{N}{u} = w \}, \]
which makes it a reflexive Banach space. Moreover, $ TH(\bou) $ is continuously embedded into $ B^{-1/2}(\bou; \mathbb{C}^3) $. In addition, the map
\[ \cross{N}{\centerdot} : TH(\bou) \longrightarrow (TH(\bou))^\ast \]
which, for $ w_1,w_2 \in TH(\bou) $, is given by
\[ \dual{\cross{N}{w_1}}{w_2} = \int_\inte \inner{(\curl u)}{\overline{v}} \, dV - \int_\inte \inner{u}{(\overline{\curl v})} \, dV, \]
where $ u,v \in \Hcurl{} $ are such that $ \cross{N}{u} = w_1 $, $ \cross{N}{v} = w_2 $, is well-defined, bounded and an isomorphism. In particular, one has that $ (\cross{N}{\centerdot})^{-1} = -\cross{N}{\centerdot} $, $ (\cross{N}{\centerdot})^\ast = -\cross{N}{\centerdot} $ and
\begin{equation}\label{for:duality_TH-TH*}
\int_\inte \inner{(\curl u)}{\overline{v}} \, dV = \int_\inte \inner{u}{(\overline{\curl v})} \, dV -\dual{\cross{N}{u}}{\cross{N}{(\cross{N}{v})}}.
\end{equation}

Before going further let us point out that, if $ w \in TH(\bou) $, $ f \in C^{0,1}(\bou) $ and $ \tilde{f} $ is any extension of $ f $ such that $ \tilde{f} \in C^{0,1}(\overline{\inte}) $, then
\[ \dual{fw}{z} = \dual{w}{\overline{f}z} = \dual{\cross{N}{(\tilde{f}u)}}{z}, \]
for any $ z \in B^{1/2}(\bou; \mathbb{C}^3) $ and any $ u \in \Hcurl{} $ such that $ \cross{N}{u} = w $.
Note that last pairing implies that $ fw \in TH(\bou) $.

Another identity to mention is the following, if $ w \in TH(\bou) $, $ f \in C^{0,1}(\bou) $ and $ \tilde{f} $ is any extension of $ f $ such that $ \tilde{f} \in C^{0,1}(\overline{\inte}) $, then
\begin{equation}\label{for:curl-product}
\dual{fw}{\cross{N}{z}} = \dual{\cross{N}{(\tilde{f}u)}}{\cross{N}{z}} = \dual{w}{\cross{N}{(\overline{f}z)}},
\end{equation}
for any $ z \in TH(\bou) $ and any $ u \in \Hcurl{} $ such that $ \cross{N}{u} = w $.

We can define the \emph{surface divergence} operator over the space $ TH(\bou) $ 
\[ \Div : TH(\bou) \longrightarrow B^{-1/2}(\bou), \]
as
\[ \Div \, w = -\inner{N}{(\curl u)}, \]
where $ u \in \Hcurl{} $ and $ \cross{N}{u} = w $. Since $ \curl u \in \Hdiv{} $, $ -\inner{N}{(\curl u)} $ makes sense and it belongs to $ B^{-1/2}(\bou) $; this operator is well-defined and bounded.

The surface divergence can be used to give an intrinsic description of the space $ TH(\bou) $. In \cite{M} Mitrea proved that there exist constants $ C_1, C_2 > 0 $ such that, for any $ w \in TH(\bou) $, the following estimates hold
\begin{gather}
\label{es:equivalent_norm_TH-C1}
\norm{w}{}{B^{-1/2}(\bou; \mathbb{C}^3)} + \norm{\Div \, w}{}{B^{-1/2}(\bou)} \leq C_1 \norm{w}{}{TH(\bou)}, \\
\label{es:equivalent_norm_TH-C2}
\norm{w}{}{TH(\bou)} \leq C_2 \left( \norm{w}{}{B^{-1/2}(\bou; \mathbb{C}^3)} + \norm{\Div \, w}{}{B^{-1/2}(\bou)} \right).
\end{gather}

\begin{lemma}\rm
There exists a positive constant $ C $ such that the following estimate holds
\begin{equation}\label{es:TH-product}
\norm{fw}{}{TH(\bou)} \leq C \norm{f}{}{C^{0,1}(\bou)} \norm{w}{}{TH(\bou)},
\end{equation}
for any $ w \in TH(\bou) $ and any $ f \in C^{0,1}(\bou) $.
\end{lemma}

\begin{proof} Consider $ \tilde{f} $ an extension of $ f $ such that $ \tilde{f} \in C^{0,1}(\overline{\inte}) $ satisfying
\[ \norm{\tilde{f}}{}{C^{0,1}(\overline{\inte})} \leq C \norm{f}{}{C^{0,1}(\bou)} \]
and any $ u \in \Hcurl{} $ such that $ \inner{N}{u} = w $.
Then
\begin{gather*}
\norm{fw}{}{TH(\bou)} \leq \norm{\tilde{f}u}{}{\Hcurl{}} \leq \norm{\tilde{f}}{}{C^{0,1}(\overline{\inte})} \norm{u}{}{\Hcurl{}} \\
\leq C \norm{f}{}{C^{0,1}(\bou)} \norm{u}{}{\Hcurl{}}.
\end{gather*}
Taking infimum in $ u \in \Hcurl{} $ one gets the estimate.
\end{proof}

\subsection{Maxwell's system as a Sch\"odinger equation}\label{sec:Me-Se}
In this section we shall transform Maxwell's equations into a Schr\"odinger-type equation. The idea of this transformation was already introduced in \cite{OS}.

There is a well known process which allows us to transform Maxwell's equations into a Schr\"odinger-type equation, in order to do so we require some extra smoothness of the coefficients, namely $ \mu, \gamma \in C^{1,1}(\overline{\inte}) $. The first step in this process  is to augment the Maxwell system with two scalar equations:
\begin{equation*}
\diver(\gamma E)=0, \qquad \diver(\mu H)=0.
\end{equation*}
The information coded in these scalar equations was already present in the initial system. In order to check this, it is enough to take divergence in each equation in (\ref{ME_without_sources}).

Next, we introduce a new system inspired in the four equations mentioned. This new system reads as
\begin{equation*}
\left\{ 
\begin{array}{l}
\diver(\gamma E) + i \omega \gamma \mu h =0\\
- \gamma^{-1} \nabla(\gamma e) + \curl E - i \omega \mu H = 0\\
\diver(\mu H) + i \omega \gamma \mu e =0\\
\mu^{-1} \nabla(\mu h) + \curl H + i \omega \gamma E = 0.\\
\end{array}
\right.
\end{equation*}
The new terms preserve the physical units of measure of the original four equations. The new system, called henceforth \emph{augmented system}, can be written as matrices in the following way
\begin{equation*}
\left[
\left(\begin{array}{c c | c c}
 &  &  & \Diver\\
 &  & D & - \Durl\\
\hline
 & \Diver &  & \\
D & \Durl &  & \\
\end{array}\right)
+\left(\begin{array}{c c | c c}
\omega \mu &  &  & D \alpha \cdot\\
 & \omega \mu I_3 & D \alpha & \\
\hline
 & D \beta \cdot & \omega \gamma & \\
D \beta & & & \omega \gamma I_3\\
\end{array}\right)
\right]\frv{h}{H}{e}{E}=0,
\end{equation*}
where $ \alpha = \mathrm{log} \, \gamma $, $ \beta = \mathrm{log} \, \mu $, $ I_j $ is $ ( j \times j ) $-identity matrix, with $ j \in \mathbb{N} $  and
\begin{gather*}
\Diver = \frac{1}{i}\left(\begin{array}{c c c}
\partial_{x^1} & \partial_{x^2} & \partial_{x^3}
\end{array}\right),\\
D = \frac{1}{i} \left( \begin{array}{c}
\partial_{x^1}\\
\partial_{x^2}\\
\partial_{x^3}
\end{array}\right),\quad
\Durl = \frac{1}{i}\left(\begin{array}{c c c}
 & -\partial_{x^3} & \partial_{x^2} \\
\partial_{x^3} &  & -\partial_{x^1} \\
-\partial_{x^2} & \partial_{x^1} & \\
\end{array}\right).
\end{gather*}
In a much more compact manner we shall express the augmented system as $ ( P + V ) X = 0 $, where
\[ P = \left(\begin{array}{c c | c c}
 &  &  & \Diver\\
 &  & D & - \Durl\\
\hline
 & \Diver &  & \\
D & \Durl &  & \\
\end{array}\right),
\qquad
V = \left(\begin{array}{c c | c c}
\omega \mu &  &  & D \alpha \cdot\\
 & \omega \mu I_3 & D \alpha & \\
\hline
 & D \beta \cdot & \omega \gamma & \\
D \beta & & & \omega \gamma I_3\\
\end{array}\right). \]

Note that $ E, H $ is a solution for Maxwell's equations, if and only if, $ X^t = \fcv{h}{H^t}{e}{E^t} $ is a solution for the augmented system and the scalar fields $ e, h $ vanish.

The next step is to rescale the augmented system, that is
\[ ( P + V ) \, 
\left( \begin{array}{c | c}
\mu^{-1/2} I_4 & \\
\hline
 & \gamma^{-1/2} I_4 \\
\end{array} \right)Y
 = \left( \begin{array}{c | c}
\gamma^{-1/2} I_4 & \\
\hline
 & \mu^{-1/2} I_4 \\
\end{array} \right) \, ( P + W )Y, \]
where
\begin{equation}\label{term:W-def}
W = \kappa I_8 + \frac{1}{2}
\left(\begin{array}{c c|c c}
 &  &  & D \alpha \cdot \\
 &  & D \alpha & D \alpha \times \\
\hline
 & D \beta \cdot &  & \\
D \beta & - D \beta \times &  & \\
\end{array}\right),
\end{equation}
with $ \kappa = \omega \mu^{1/2} \gamma^{1/2}$. We shall call
\[ ( P + W ) Y = 0 \]
the \emph{rescaled system}.

The advantage of rescaling is that
\begin{align}
0 &= ( P + W )( P - W^t ) Z = (-\Delta I_8 + Q) Z, \label{eq:schr1}\\
0 &= ( P - W^t )( P + W ) Z' = (-\Delta I_8 + Q') Z', \label{eq:schr2}\\
0 &= ( P + W^* )( P - \overline{W} ) \hat Z = (-\Delta I_8 + \hat Q) \hat Z, \label{eq:schr3}
\end{align}
where $ Q, Q', \hat Q $ are zeroth-order terms. Here $ W^t $ denotes the transposed of $ W $ and $ W^\ast $ stands for $ \overline{W^t} $. No first order terms appear in (\ref{eq:schr1}), (\ref{eq:schr2}) and (\ref{eq:schr3}), giving as a result a Schr\"odinger-type equation.
Mind
\begin{equation}
Q = -PW^t + WP - WW^t\label{def:matrixQ-0th}.
\end{equation}
Note that if $ Z $ is a solution for (\ref{eq:schr1}) in $ \inte $, then $ Y = ( P - W^t ) Z $ is a solution for the rescaled system in $ \inte $, hence
\[ X = \left( \begin{array}{c | c}
\mu^{-1/2} I_4 &  \\
\hline
 & \gamma^{-1/2} I_4 \\
\end{array} \right) Y \]
is a solution for the augmented system. In the same manner, if $ \hat Z $ is a solution for (\ref{eq:schr3}), then $ \hat Y = (P - \overline{W})\hat Z $ is a solution for $ (P + W^*) \hat Y = 0 $ in $ \inte $.

For later uses,
\begin{gather}
Q = \frac{1}{2} \fbM{\subM{\Delta \alpha}{}{}{2 \nabla^2 \alpha - \Delta \alpha I_3}}{}
{}{\subM{\Delta \beta}{}{}{2 \nabla^2 \beta - \Delta \beta I_3}}\nonumber \\
\label{ter:matrixQ-0th}
- \fbM{(\kappa^2 + \frac{1}{4} (\inner{D\alpha}{D\alpha})) I_4}{\subM{}{2 D \kappa \cdot}{2 D \kappa}{}}
{\subM{}{2 D \kappa \cdot}{2 D \kappa}{}}{(\kappa^2 + \frac{1}{4} (\inner{D\beta}{D\beta})) I_4},
\end{gather}
\begin{gather}
Q' = - \frac{1}{2} \fbM{\subM{\Delta \beta}{}{}{2 \nabla^2 \beta - \Delta \beta I_3}}{}
{}{\subM{\Delta \alpha}{}{}{2 \nabla^2 \alpha - \Delta \alpha I_3}}\nonumber \\
\label{ter:matrixQ'-0th}
- \fbM{(\kappa^2 + \frac{1}{4} (\inner{D\beta}{D\beta})) I_4}{\subM{}{}{}{2D \kappa \times}}
{\subM{}{}{}{-2 D \kappa \times}}{(\kappa^2 + \frac{1}{4} (\inner{D\alpha}{D\alpha})) I_4}
\end{gather}
and
\begin{gather}
\hat Q = \frac{1}{2} \fbM{\subM{- \Delta \beta}{}{}{- 2 \nabla^2 \beta + \Delta \beta I_3}}{}
{}{\subM{- \Delta \overline{\alpha}}{}{}{- 2 \nabla^2 \overline{\alpha} + \Delta \overline{\alpha} I_3}}\nonumber \\
\label{ter:matrixQ^-0th}
- \fbM{(\overline{\kappa}^2 - \frac{1}{4} (\inner{D\beta}{D\beta})) I_4}{\subM{}{}{}{-2 D \overline{\kappa} \times}}
{\subM{}{}{}{2 D \overline{\kappa} \times}}{(\overline{\kappa}^2 - \frac{1}{4} (\inner{D \overline{\alpha}}{D \overline{\alpha}})) I_4}
\end{gather}
with $\nabla^2 f=(\partial^2_{x_j, x_k} f)^3_{j,k=1}$.

The computations needed to get (\ref{ter:matrixQ-0th}) and (\ref{ter:matrixQ'-0th}) can be found in \cite{COSa}. The same kind of computations gives (\ref{ter:matrixQ^-0th}).

In order to make as concise as possible the presentation of our proof, we introduce some additional notation. Let $ Y, Z $ be in the form
\[ Y = \fcv{f^1}{(u^1)^t}{f^2}{(u^2)^t}^t, \qquad Z = \fcv{g^1}{(v^1)^t}{g^2}{(v^2)^t}^t, \]
define
\begin{equation*}
\Inner{Y}{Z}_\inte = \sum_{j=1}^{2} \left( \int_{\inte} f^j \overline{g^j} \, dV + \int_{\inte} \inner{u^j}{\overline{v^j}} \, dV \right),
\end{equation*}
\begin{equation*}
\Inner{Y}{Z}_{\bou} =\sum_{j=1}^{2} \left( \int_{\bou} f^j \overline{g^j} \, dA + \int_{\bou} \inner{u^j}{\overline{v^j}} \, dA \right).
\end{equation*}
In the first identity we are assuming $ f^j, g^j \in C^\infty(\overline{\inte}) $ and $ u^j, v^j \in \mathcal{X}\mathbb{E}|_{\overline{\inte}} $ with $ j = 1, 2 $, while in the second identity $ f^j, g^j \in C^\infty(\bou) $ and $ u^j, v^j \in \mathcal{X}\mathbb{E}|_{\bou} $ with $ j = 1, 2 $.
The following integration by parts holds
\[ \Inner{PY}{Z}_\inte = \Inner{P_N Y}{Z|_{\bou}}_{\bou} + \Inner{Y}{PZ}_\inte. \]
Here, when $ A $ is a (possibly complex) vector field we denote
\begin{equation}\label{for:P_A}
P_A = \frac{1}{i} \left( \begin{array}{c c|c c}
 &  &  & A \cdot \\
 &  & A & - A \times \\
\hline
 & A \cdot &  & \\
A & A \times &  & \\
\end{array}\right).
\end{equation}

Finally, for elements $ Y $ in the form given above we define, for $ |s| > 0 $,
\[ \norm{Y}{}{H^s(\inte; \mathcal{Y})} = \sum_{j = 1, 2} \left( \norm{f^j}{}{H^s(\inte)} + \norm{u^j}{}{H^s(\inte; \mathbb{C}^3)} \right), \]
and
\[ \norm{Y}{}{L^2(\inte; \mathcal{Y})} = \sum_{j = 1, 2} \left( \norm{f^j}{}{L^2(\inte)} + \norm{u^j}{}{L^2(\inte; \mathbb{C}^3)} \right). \]
Here and henceforth, $ \mathcal{Y} $ denotes the space $ \mathbb{C} \times \mathbb{C}^3 \times \mathbb{C} \times \mathbb{C}^3 $. On the other hand, we define, for $ 0 < |s| < 1 $,
\[ \norm{Y}{}{B^s(\bou; \mathcal{Y})} = \sum_{j = 1, 2} \left( \norm{f^j}{}{B^s(\bou)} + \norm{u^j}{}{B^s(\bou; \mathbb{C}^3)} \right), \]
and
\[ \norm{Y}{}{L^2(\bou; \mathcal{Y})} = \sum_{j = 1, 2} \left( \norm{f^j}{}{L^2(\bou)} + \norm{u^j}{}{L^2(\bou; \mathbb{C}^3)} \right). \]

\section{A log-type estimate}\label{sec:log-type_estimate}
In this section we prove Theorem \ref{th:stability}. The proof consists of three ingredients, namely, a suitable estimate relating the electromagnetic coefficients in $ \inte $ with their corresponding Cauchy data sets, the constructions of special solutions and the use of a Carleman estimate.

\begin{lemma}\label{le:1-suitable-estimate}\rm
Let $ \mu_j, \gamma_j $ belong to $ C^{0,1}(\overline{\inte}) $. Then one has that, for any $ Y_1 $ given by
\[ Y_1 = \fcv{0}{\mu_1^{1/2} H_1^ t}{0}{\gamma_1^{1/2} E_1^t}^t \]
with $ E_1, H_1 \in \Hcurl{} $ solution for (\ref{ME_without_sources}) in $ \inte $ with coefficients $ \mu_1, \gamma_1 $, and any
\[ Y_2 = \fcv{f^1}{(u^1)^t}{f^2}{(u^2)^t}^t \in H^1(\inte) \times \Hcurl{} \times H^1(\inte) \times \Hcurl{} \]
solution for $ ( P + W^\ast_2 ) Y_2 = 0 $ in $ \inte $; the following estimate holds:
\begin{gather*}
|\Inner{Y_1}{PY_2}_\inte - \Inner{PY_1}{Y_2}_\inte| \\
\leq C \delta_C(C_1,C_2) \left( \norm{\mu_2^{-1/2}}{}{C^{0,1}(\bou)} \norm{g_2}{}{B^{1/2}(\bou)} + \norm{\gamma_2^{1/2}}{}{C^{0,1}(\bou)} \norm{z_1}{}{TH(\bou)} \right. \\
\left. + \norm{\gamma_2^{-1/2}}{}{C^{0,1}(\bou)} \norm{g_1}{}{B^{1/2}(\bou)} + \norm{\mu_2^{1/2}}{}{C^{0,1}(\bou)} \norm{z_2}{}{TH(\bou)} \right) \norm{\cross{N}{E_1}}{}{TH(\bou)} \\
+ C \left( \norm{\cross{N}{E_1}}{}{TH(\bou)} + \norm{\cross{N}{H_1}}{}{TH(\bou)} \right) \\
\times\left( \norm{\mu_1^{-1/2} - \mu_2^{-1/2}}{}{C^{0,1}(\bou)} \norm{g_2}{}{B^{1/2}(\bou)} + \norm{\mu_1^{1/2} - \mu_2^{1/2}}{}{C^{0,1}(\bou)} \norm{z_2}{}{TH(\bou)} \right. \\
\left. + \norm{\gamma_1^{1/2} - \gamma_2^{1/2}}{}{C^{0,1}(\bou)} \norm{z_1}{}{TH(\bou)} + \norm{\gamma_1^{-1/2} - \gamma_2^{-1/2}}{}{C^{0,1}(\bou)} \norm{g_1}{}{B^{1/2}(\bou)} \right).
\end{gather*}
Here $ g_1, g_2 \in B^{1/2}(\bou) $ stand for $ g_1 = f^1|_{\bou}, \, g_2 = f^2|_{\bou} $ while $ z_1, z_2 \in TH(\bou) $ stand for $ z_1 = \cross{N}{u^1}, \, z_2 = \cross{N}{u^2} $. Here $ C_j$, with $ j=1,2 $, stands for the Cauchy data set corresponding to $ \mu_j, \gamma_j $. Recall that $ W_2 $ is the matrix (\ref{term:W-def}) associated to $ \mu_2, \gamma_2 $.
\end{lemma}

Along these notes, $ \mu_j $ and $ \gamma_j $ should be understood either as themselves or as their traces, according to the context.

\begin{proof}
Let $ L $ be
\[ L = \fcv{0}{\mu_2^{1/2} H_2^t}{0}{\gamma_2^{1/2} E_2^t}^t, \]
with $ E_2, H_2 \in \Hcurl{} $ an arbitrary solution for (\ref{ME_without_sources}) with coefficients $ \mu_2, \gamma_2 $. Since $ ( P + W_2^\ast ) Y_2 = 0 $ and $ ( P + W_2 ) L = 0 $, one has that $ \Inner{L}{PY_2}_\inte = \Inner{PL}{Y_2}_\inte $, hence
\[ \Inner{Y_1}{PY_2}_\inte - \Inner{PY_1}{Y_2}_\inte = \Inner{Y_1 - L}{PY_2}_\inte - \Inner{P(Y_1 - L)}{Y_2}_\inte. \]
On the other hand, we have, using (\ref{def:normal-component}), (\ref{for:div-product}), (\ref{for:duality_TH-TH*}) and (\ref{for:curl-product}), that
\begin{gather*}
\Inner{Y_1 - L}{PY_2}_\inte - \Inner{P(Y_1 - L)}{Y_2}_\inte =\\
= i\dual{\inner{N}{(\mu_1 H_1 - \mu_2 H_2)}}{\mu_2^{-1/2} g_2} + i\dual{\inner{N}{(\mu_1 H_1)}}{(\mu_1^{-1/2} - \mu_2^{-1/2}) g_2}\\
+ i\dual{\inner{N}{(\gamma_1 E_1 - \gamma_2 E_2)}}{\overline{\gamma_2^{-1/2}} g_1} + i\dual{\inner{N}{(\gamma_1 E_1)}}{(\overline{\gamma_1^{-1/2} - \gamma_2^{-1/2}}) g_1}\\
- i\dual{\cross{N}{(H_1-H_2)}}{\cross{N}{(\mu_2^{1/2}z_2)}} - i\dual{\cross{N}{H_1}}{\cross{N}{((\mu_1^{1/2} - \mu_2^{1/2})z_2)}}\\
+ i\dual{\cross{N}{(E_1-E_2)}}{\cross{N}{(\overline{\gamma_2^{1/2}}z_1)}} + i\dual{\cross{N}{E_1}}{\cross{N}{((\overline{\gamma_1^{1/2} - \gamma_2^{1/2}})z_1)}}.
\end{gather*}
Furthermore, from Maxwell's equations one deduces that
\[ \inner{N}{(\gamma_j E_j)} = \frac{1}{i\omega} \Div \, (\cross{N}{H_j}), \quad \inner{N}{(\mu_j H_j)} = -\frac{1}{i\omega} \Div \, (\cross{N}{E_j}), \]
for $ j = 1,2 $. Hence, denoting $ \cross{N}{E_j} = T_j $ and $ \cross{N}{H_j} = S_j $ we obtain
\begin{gather*}
\Inner{Y_1}{PY_2}_\inte - \Inner{PY_1}{Y_2}_\inte =\\
= -\frac{1}{\omega} \dual{\Div (T_1 - T_2)}{\mu_2^{-1/2} g_2} -\frac{1}{\omega} \dual{\Div\, T_1}{(\mu_1^{-1/2} - \mu_2^{-1/2}) g_2}\\
+ \frac{1}{\omega} \dual{\Div(S_1 - S_2)}{\overline{\gamma_2^{-1/2}} g_1} + \frac{1}{\omega} \dual{\Div\, S_1}{(\overline{\gamma_1^{-1/2} - \gamma_2^{-1/2}}) g_1}\\
- i\dual{S_1 - S_2}{\cross{N}{(\mu_2^{1/2}z_2)}} - i\dual{S_1}{\cross{N}{((\mu_1^{1/2} - \mu_2^{1/2})z_2)}}\\
+ i\dual{T_1 - T_2}{\cross{N}{(\overline{\gamma_2^{1/2}}z_1)}} + i\dual{T_1}{\cross{N}{((\overline{\gamma_1^{1/2} - \gamma_2^{1/2}})z_1)}}.
\end{gather*}
By using the appropriate dualities and the estimates (\ref{es:besov_product}), (\ref{es:TH-product}) and (\ref{es:equivalent_norm_TH-C1}) we get
\begin{gather*}
|\Inner{Y_1}{PY_2}_\inte - \Inner{PY_1}{Y_2}_\inte| \leq\\ 
\leq C \left( \norm{T_1 - T_2}{}{TH(\bou)} + \norm{S_1 - S_2}{}{TH(\bou)} \right) \left( \norm{\mu_2^{-1/2}}{}{C^{0,1}(\bou)} \norm{g_2}{}{B^{1/2}(\bou)} \right. \\
+ \norm{\gamma_2^{1/2}}{}{C^{0,1}(\bou)} \norm{z_1}{}{TH(\bou)} + \norm{\gamma_2^{-1/2}}{}{C^{0,1}(\bou)} \norm{g_1}{}{B^{1/2}(\bou)}\\
\left. + \norm{\mu_2^{1/2}}{}{C^{0,1}(\bou)} \norm{z_2}{}{TH(\bou)} \right) + C \left( \norm{\mu_1^{-1/2} - \mu_2^{-1/2}}{}{C^{0,1}(\bou)} \norm{g_2}{}{B^{1/2}(\bou)} \right. \\
+ \norm{\gamma_1^{1/2} - \gamma_2^{1/2}}{}{C^{0,1}(\bou)} \norm{z_1}{}{TH(\bou)} + \norm{\gamma_1^{-1/2} - \gamma_2^{-1/2}}{}{C^{0,1}(\bou)} \norm{g_1}{}{B^{1/2}(\bou)}\\
\left. + \norm{\mu_1^{1/2} - \mu_2^{1/2}}{}{C^{0,1}(\bou)} \norm{z_2}{}{TH(\bou)} \right)\left( \norm{T_1}{}{TH(\bou)} + \norm{S_1}{}{TH(\bou)} \right).
\end{gather*}
This estimate holds for all $ (T_2, S_2) \in C_2 $, since $ E_2, H_2 $ was chosen to be an arbitrary solution for (\ref{ME_without_sources}) with coefficients $ \mu_2, \gamma_2 $. Finally, the wanted estimate is a consequence of Definition \ref{def:Cdistance}.
\end{proof}

\begin{proposition}\label{le:2-suitable-estimate}\rm
Let $ \gamma_1, \mu_1 $ and $ \gamma_2, \mu_2 $ be in the class $ B $-stable coefficients on the boundary at frequency $ \omega $, with $ B $ as in Theorem \ref{th:stability}. Then, there exists a constant $ C(M) $ such that, for any $ Z_1 \in H^1(\inte; \mathcal{Y}) $ satisfying $ Y_1 = (P - W_1^t) Z_1 $ with $ Y_1 $ as in Lemma \ref{le:1-suitable-estimate} and any $ Y_2 \in H^1(\inte; \mathcal{Y}) $ as in Lemma \ref{le:1-suitable-estimate}, one has
\begin{gather}
|\Inner{(Q_1 - Q_2)Z_1}{Y_2}_\inte| \leq C\, B\big( \delta_C(C_1, C_2) \big) \norm{Z_1}{}{H^1(\inte; \mathcal{Y})} \norm{Y_2}{}{H^1(\inte; \mathcal{Y})}\nonumber\\
+\, C \, B\big( \delta_C(C_1, C_2) \big) \left( \norm{E_1}{}{\Hcurl{}} + \norm{H_1}{}{\Hcurl{}} \right)\nonumber\\
\label{es:suitablest}
\times \left( \norm{f^1}{}{H^1(\inte)} + \norm{u^1}{}{\Hcurl{}} + \norm{u^2}{}{\Hcurl{}} + \norm{f^2}{}{H^1(\inte)} \right).
\end{gather}
Here $ Q_j $ is the matrix (\ref{def:matrixQ-0th}) associated to $ \mu_j, \gamma_j $ with $ j=1,2 $.\nopagebreak[4]
\end{proposition}

\begin{proof}
From (\ref{def:matrixQ-0th}) one has
\begin{gather*}
\Inner{(Q_1 - Q_2)Z_1}{Y_2}_\inte = - \Inner{P(W^t_1 - W^t_2)Z_1}{Y_2}_\inte\\
+ \Inner{(W_1 - W_2)PZ_1}{Y_2}_\inte - \Inner{(W_1W^t_1 - W_2W^t_2)Z_1}{Y_2}_\inte\\
= \Inner{(W_1^t - W_2^t)Z_1}{P_N Y_2}_{\bou} - \Inner{W^t_1 Z_1}{PY_2}_\inte + \Inner{W^t_2 Z_1}{PY_2}_\inte\\
  + \Inner{W_1(P - W_1^t) Z_1}{Y_2}_\inte - \Inner{PZ_1}{W_2^\ast Y_2}_\inte + \Inner{W^t_2 Z_1}{W_2^\ast Y_2}_\inte\\
= \Inner{(W_1^t - W_2^t)Z_1}{P_N Y_2}_{\bou} + \Inner{(P-W_1^t)Z_1}{PY_2}_\inte  + \Inner{W_1(P-W_1^t)Z_1}{Y_2}_\inte\\
= \Inner{(W_1^t - W_2^t)Z_1}{P_N Y_2}_{\bou} + \Inner{Y_1}{PY_2}_\inte - \Inner{PY_1}{Y_2}_\inte.
\end{gather*}
In order to get the penultimate identity, we used twice that $ (P + W_2^\ast)Y_2 = 0 $. In the last one, we used that $ Y_1 = (P-W^t_1)Z_1 $ and that $ (P + W_1)Y_1 = 0 $.

It is a straight forward computation to check the next estimate
\begin{gather*}
|\Inner{(W_1^t - W_2^t)Z_1}{P_N Y_2}_{\bou}| \leq C \left( \norm{\kappa_1 - \kappa_2}{}{L^\infty(\bou)} + \norm{\nabla(\beta_1 - \beta_2)}{}{L^\infty(\bou; \mathbb{C}^3)} \right.\\
\left. + \norm{\nabla(\alpha_1 - \alpha_2)}{}{L^\infty(\bou; \mathbb{C}^3)} \right) \norm{Z_1}{}{L^2(\bou; \mathcal{Y})} \norm{Y_2}{}{L^2(\bou; \mathcal{Y})}.
\end{gather*}
Here, as usually, the norm of $ L^\infty (\bou; \mathbb{C}^3) $ is
\[ \norm{w}{2}{L^\infty(\bou; \mathbb{C}^3)} = \sum_{j=1}^3 \norm{w^{(j)}}{2}{L^\infty(\bou)}, \]
for any $ w \in \mathcal{X}(\mathbb{E})|_{\bou} $.

It is a routine computation to check that, on one hand
\begin{gather*}
\begin{align*}
\norm{\kappa_1 - \kappa_2}{}{L^\infty(\bou)} &\leq C \left( \norm{\gamma_2 - \gamma_1}{}{L^\infty(\bou)} + \norm{\mu_2 - \mu_1}{}{L^\infty(\bou)} \right)\\
&\leq C\, B\big( \delta_C(C_1, C_2) \big),
\end{align*}\\
\begin{align*}
\norm{\nabla(\alpha_1 - \alpha_2)}{}{L^\infty(\bou; \mathbb{C}^3)} &\leq C \left( \norm{\gamma_1 - \gamma_2}{}{L^\infty(\bou)} + \norm{\nabla(\gamma_1 - \gamma_2)}{}{L^\infty(\bou; \mathbb{C}^3)} \right)\\
&\leq C\, B\big( \delta_C(C_1, C_2) \big),
\end{align*}\\
\begin{align*}
\norm{\nabla(\beta_1 - \beta_2)}{}{L^\infty(\bou; \mathbb{C}^3)} &\leq C \left( \norm{\mu_1 - \mu_2}{}{L^\infty(\bou)} + \norm{\nabla(\mu_1 - \mu_2)}{}{L^\infty(\bou; \mathbb{C}^3)} \right)\\
&\leq C\, B\big( \delta_C(C_1, C_2) \big).
\end{align*}
\end{gather*}
and on the other hand,
\begin{gather*}
\norm{\mu_2^{-1/2}}{}{C^{0,1}(\bou)} + \norm{\gamma_2^{-1/2}}{}{C^{0,1}(\bou)} + \norm{\mu_2^{1/2}}{}{C^{0,1}(\bou)} + \norm{\gamma_2^{1/2}}{}{C^{0,1}(\bou)} \leq C\\
\begin{align*}
\norm{\mu_1^{-1/2} - \mu_2^{-1/2}}{}{C^{0,1}(\bou)} &\leq C \norm{\mu_1 - \mu_2}{}{C^{0,1}(\bou)} \leq C\, B\big( \delta_C(C_1, C_2) \big),\\
\norm{\gamma_1^{-1/2} - \gamma_2^{-1/2}}{}{C^{0,1}(\bou)} &\leq C \norm{\gamma_1 - \gamma_2}{}{C^{0,1}(\bou)} \leq C\, B\big( \delta_C(C_1, C_2) \big),\\
\norm{\mu_1^{1/2} - \mu_2^{1/2}}{}{C^{0,1}(\bou)} &\leq C \norm{\mu_1 - \mu_2}{}{C^{0,1}(\bou)} \leq C\, B\big( \delta_C(C_1, C_2) \big),\\
\norm{\gamma_1^{1/2} - \gamma_2^{1/2}}{}{C^{0,1}(\bou)} &\leq C \norm{\gamma_1 - \gamma_2}{}{C^{0,1}(\bou)} \leq C\, B\big( \delta_C(C_1, C_2) \big).
\end{align*}
\end{gather*}
Putting together all these estimates and Lemma \ref{le:1-suitable-estimate}, we get
\begin{gather*}
|\Inner{(Q_1 - Q_2)Z_1}{Y_2}_\inte| \leq C\, B\big( \delta_C(C_1, C_2) \big) \norm{Z_1}{}{B^{1/2}(\bou; \mathcal{Y})} \norm{Y_2}{}{B^{1/2}(\bou; \mathcal{Y})}\\
+\, C \, B\big( \delta_C(C_1, C_2) \big) \left( \norm{\cross{N}{E_1}}{}{TH(\bou)} + \norm{\cross{N}{H_1}}{}{TH(\bou)} \right)\\
\times \left( \norm{g_1}{}{B^{1/2}(\bou)} + \norm{z_1}{}{TH(\bou)} + \norm{g_2}{}{B^{1/2}(\bou)} + \norm{z_2}{}{TH(\bou)} \right),
\end{gather*}
hence we deduce the estimate given in the statement.
\end{proof}

\subsection{Construction of special solutions}
Here we construct two kinds of special solutions, one for the Schr\"odinger-type equation and another one for $ (P + W^*) Y = 0 $. The first one was already constructed in \cite{OS} but we give here the proof in order to keep track the constants. The second kind of solution is inspired on the solutions given in \cite{KSaU}.

Let $ B(O; \rho) $ be the open ball centered at the origin $ O $ with radius $ \rho > 0 $ and such that $ \overline{\inte} \subset B(O; \rho) $. Sometimes $ B(O; \rho) $ will be denoted by $ B_\rho $ to simplify the notation. Let $ \varepsilon_0 $ and $ \mu_0 $ denote the electric and magnetic constants, respectively. Extend the coefficients $ \gamma, \mu $ defined in $ \inte $ to functions in $ \mathbb{E} $ --still denoted by $ \gamma, \mu $--, preserving their smoothness and in such a way that $ \gamma-\varepsilon_0, \mu-\mu_0 $ have compact support in $ \overline{B(O;\rho)} $ (regarding to extensions see \cite{St}). Note two simple facts. Firstly, the extensions still satisfy the a priori bound and the a priori ellipticity condition in $ \mathbb{E} $. Secondly, the extensions of the matrices (\ref{ter:matrixQ-0th}), (\ref{ter:matrixQ'-0th}) (\ref{ter:matrixQ^-0th}) --still denoted by $ Q, Q', \hat Q $-- satisfy that $ \omega^2 \varepsilon_0 \mu_0 I_8 + Q $, $ \omega^2 \varepsilon_0 \mu_0 I_8 + Q' $ and $ \omega^2 \varepsilon_0 \mu_0 I_8 + \hat Q $ have compact support in $ \overline{B(O;\rho)} $.

We shall construct solutions for (\ref{eq:schr1}) in $ \mathbb{E} $ with the form of a complex geometrical optic solution (CGO solution for short), that is, in the form
\[ Z = e^{i \zeta \cdot x} (L + R), \]
with $ L = L(\zeta) $ constant and $ \zeta \in \mathbb{C}^3 $.

\begin{lemma}\label{le:CGO-sch} \rm
Let $ \delta $ be a constant such that $ -1 < \delta < 0 $ and let $ \zeta \in \mathbb{C}^3 $ be such that $ \zeta \cdot \zeta = \omega^2 \varepsilon_0 \mu_0 $ with
\begin{equation*}
|\zeta| > C \sum_{j,k = 1}^8 \norm{(\omega^2 \varepsilon_0 \mu_0 I_8 + Q)_j^k}{}{L^\infty(B_\rho)}
\end{equation*}
for some constant $ C = C(\delta, \rho) > 0 $.
Then, there exists a
\[ Z = e^{i \zeta \cdot x} (L + R) \]
solution for $ (-\Delta I_8 + Q) Z = 0 $ in $ \mathbb{E} $, with $ Z|_\inte \in H^2(\inte; \mathcal{Y}) $, $ L = L(\zeta) $ constant and $ R = R(\zeta) $ satisfying 
\begin{align}
\label{es:decay_R}
\norm{R}{}{L^2_\delta \mathcal{Y}} &\leq \frac{C(\delta, \rho)}{|\zeta|} |L| \sum_{j,k=1}^8 \norm{(\omega^2 \varepsilon_0 \mu_0 I_8 + Q)_j^k}{}{L^\infty(B_\rho)},\\
\label{es:decay_PR}
\norm{PR}{}{L^2_\delta \mathcal{Y}} &\leq C(\delta, \rho) \sum_{j,k=1}^8 \norm{(\omega^2 \varepsilon_0 \mu_0 I_8 + Q)_j^k}{}{L^\infty(B\rho)} \big(|L| + \norm{R}{}{L^2_\delta \mathcal{Y}} \big).
\end{align}
\end{lemma}

The norm in the lemma is
\[ \norm{Y}{}{L^2_\lambda \mathcal{Y}} = \sum_{j = 1, 2} \left( \norm{f^j}{}{L^2_\lambda} + \norm{u^j}{}{L^2_\lambda(\mathbb{E}; \mathbb{C}^3)} \right), \]
with $ Y = \fcv{f^1}{(u^1)^t}{f^2}{(u^2)^t}^t $ and
\[ \norm{f}{2}{L^2_\lambda} = \int_{\mathbb{R}^3} (1 + |x|^2)^\lambda |f|^2 \, dx, \]
for $ 0< |\lambda| < 1 $.

\begin{proof}
It is an easy matter to check that
\[ e^{-i \zeta \cdot x} (-\Delta I_8) e^{i \zeta \cdot x} (L + R) = [-\Delta -2i\zeta \cdot \nabla + \zeta \cdot \zeta]I_8(L + R), \]
hence, if $ Z = e^{i \zeta \cdot x} (L + R) $ is a solution for (\ref{eq:schr1}), then $ R $ solves
\begin{equation}\label{eq:ansatz-CGO}
((-\Delta -2i\zeta \cdot \nabla + \omega^2 \varepsilon_0 \mu_0)I_8 + Q) R = -(\omega^2 \varepsilon_0 \mu_0 I_8 + Q) L.
\end{equation}
We will use this equation as the starting point for the construction of the CGO solutions.

It is a very well known fact that the operator $ G_\zeta $ --standing for the convolution with the fundamental solution for Fadeev's operator $ (-\Delta -2i\zeta \cdot \nabla) $-- with $ |\zeta| > 1 $ satisfy the estimate
\begin{align}\label{es:syl-uhl}
\norm{G_\zeta f}{}{L^2_{\delta}} &\leq \frac{C( \delta)}{|\zeta|} \norm{f}{}{L^2_{\delta + 1}},\\
\label{es:brown}
\norm{\partial_{x^j} G_\zeta f}{}{L^2_{\delta}} &\leq C( \delta) \norm{f}{}{L^2_{\delta + 1}},
\end{align}
for any $ f \in L^2_{\delta + 1} $ with $ j = 1, 2, 3 $. Estimate (\ref{es:syl-uhl}) was proven in \cite{SyU}, while (\ref{es:brown}) was proven in \cite{B}.

Denote $ F_\zeta = G_\zeta I_8 $. Applying $ F_\zeta $ to both sides of (\ref{eq:ansatz-CGO}) we get
\begin{equation}\label{eq:post-ansatz-CGO}
(I_8 + F_\zeta(\omega^2 \varepsilon_0 \mu_0 I_8 + Q)) R = - F_\zeta(\omega^2 \varepsilon_0 \mu_0 I_8 + Q) L.
\end{equation}

On the other hand, we can estimate
\begin{gather*}
\norm{F_\zeta ( \omega^2 \varepsilon_0 \mu_0 I_8 + Q ) R}{}{L^2_\delta \mathcal{Y}} \leq \frac{C(\delta, \rho)}{|\zeta|} \sum_{j,k=1}^8 \norm{(\omega^2 \varepsilon_0 \mu_0 I_8 + Q)_j^k}{}{L^\infty(B_\rho)} \norm{R}{}{L^2_{\delta}\mathcal{Y}},
\end{gather*}
where we applied (\ref{es:syl-uhl}) and used the fact that $ \omega^2 \varepsilon_0 \mu_0 I_8 + Q $ has compact support in $ \overline{B(O; \rho)} $.

Then, assuming
\begin{equation*}
|\zeta| > C(\delta, \rho) \sum_{j,k = 1}^8 \norm{(\omega^2 \varepsilon_0 \mu_0 I_8 + Q)_j^k}{}{L^\infty(B_\rho)},
\end{equation*}
the operator $ ( I_8 + F_\zeta(\omega^2 \varepsilon_0 \mu_0 I_8 + Q) )^{-1} $ is bounded in $ L^2_\delta \mathcal{Y} $ and
\[ R = - (I_8 + F_\zeta(\omega^2 \varepsilon_0 \mu_0 I_8 + Q))^{-1} F_\zeta(\omega^2 \varepsilon_0 \mu_0I_8 + Q) L, \]
with
\begin{equation*}
\norm{R}{}{L^2_\delta \mathcal{Y}} \leq \frac{C(\delta, \rho)}{|\zeta|} |L| \sum_{j,k=1}^8 \norm{(\omega^2 \varepsilon_0 \mu_0 I_8 + Q)_j^k}{}{L^\infty(B_\rho)}.
\end{equation*}
Here we used again (\ref{es:syl-uhl}) and the fact that $ \omega^2 \varepsilon_0 \mu_0 I_8 + Q $ has compact support in $ \overline{B(O; \rho)} $. This compactness is crucial in our arguments.

On the other hand, from (\ref{eq:post-ansatz-CGO}) and (\ref{es:brown}) we deduce that
\begin{align*}
\norm{PR}{}{L^2_\delta \mathcal{Y}} &\leq \norm{P [F_\zeta(\omega^2 \varepsilon_0 \mu_0 I_8 + Q)(R + L)]}{}{L^2_\delta\mathcal{Y}} \\
&\leq C(\delta, \rho) \sum_{j,k=1}^8 \norm{(\omega^2 \varepsilon_0 \mu_0 I_8 + Q)_j^k}{}{L^\infty(B_\rho)} \big(|L| + \norm{R}{}{L^2_\delta \mathcal{Y}} \big).
\end{align*}
\end{proof}

The arguments given in the above proof were used in \cite{SyU} for scalar equations of the same type as (\ref{eq:ansatz-CGO}). We exploit below the uniqueness for these scalar equations.

We will use the CGO solutions constructed for (\ref{eq:schr1}) to produce analogous solutions for Maxwell's equations. The procedure follows the ideas exposed in Section \ref{sec:Me-Se}, using the decoupled scalar equations in (\ref{eq:schr2}).

\begin{proposition}\label{pro:CGO-sch} \rm
Let $ \delta $ be a constant such that $ -1 < \delta < 0 $ and let $ \zeta \in \mathbb{C}^3 $ be such that $ \zeta \cdot \zeta = \omega^2 \varepsilon_0 \mu_0 $ with
\[ |\zeta| > C(\delta, \rho)\left(\sum_{j = 1,2} \norm{\omega^2 \varepsilon_0 \mu_0 + q_j}{}{L^\infty(B_\rho)} + \sum_{j,k=1}^8 \norm{(\omega^2 \varepsilon_0 \mu_0 I_8 + Q)_j^k}{}{L^\infty(B_\rho)} \right), \]
where
\begin{equation*}
q_1 = -\frac{1}{2} \Delta \beta - \kappa^2 - \frac{1}{4} (\inner{D\beta}{D\beta}), \quad q_2 = - \frac{1}{2} \Delta \alpha - \kappa^2 - \frac{1}{4} (\inner{D\alpha}{D\alpha}).
\end{equation*}
Then, there exists a
\[ Z = e^{i \zeta \cdot x} (L + R) \]
solution of $ (-\Delta I_8 + Q) Z = 0 $ in $ \mathbb{E} $, with $ Z|_\inte \in H^2(\inte; \mathcal{Y}) $,
\[ L = \frac{1}{|\zeta|} \frv{\zeta \cdot a}{\omega \varepsilon_0^{1/2} \mu_0^{1/2} b}{\zeta \cdot b}{\omega \varepsilon_0^{1/2} \mu_0^{1/2} a}, \]
for $ a, b $ constant complex vector fields, and $ R $ satisfying
\begin{align*}
\norm{R}{}{L^2_\delta \mathcal{Y}} &\leq \frac{C(\delta, \rho)}{|\zeta|} |L| \sum_{j,k=1}^8 \norm{(\omega^2 \varepsilon_0 \mu_0 I_8 + Q)_j^k}{}{L^\infty(B_\rho)}.
\end{align*}
Furthermore, $ Y = (P - W^t)Z $ is solution for $ (P + W)Y = 0 $ in $ \mathbb{E} $ and it reads
\[ Y = \fcv{0}{\mu^{1/2} H^t}{0}{\gamma^{1/2} E^t}^t \]
with $ E, H $ solution for (\ref{ME_without_sources}) in $ \mathbb{E} $.
\end{proposition}

\begin{proof}
Let $ Y $ be defined by $ Y = (P - W^t)Z $, with $ Z $ the solution constructed in Lemma \ref{le:CGO-sch}. If we denote
\[ Y = \fcv{f^1}{(u^1)^t}{f^2}{(u^2)^t}^t, \]
we will prove that $ f^1 = f^2 = 0 $.

Note that $ Y $ solves (\ref{eq:schr2}) weakly, with $ f^j $ solving the following decoupled equation:
\begin{equation*}
(- \Delta + q_j)f^j = 0, \qquad j=1, 2
\end{equation*}
in the weak sense.

Denoting $ L, R $ from Lemma \ref{le:CGO-sch} as
\begin{equation}\label{ter:Lform_Rform}
L = \fcv{l^1}{(L^1)^t}{l^2}{(L^2)^t}^t, \quad R = \fcv{r^1}{(R^1)^t}{r^2}{(R^2)^t}^t;
\end{equation}
the functions $ f^j $ can be expressed as $ f^j = e^{i\zeta \cdot x} (m_j + s_j) $
with
\begin{align*}
m_1 &= \zeta \cdot L^2 - \kappa l^1 \quad & s_1 = -\frac{1}{2} D\beta \cdot L^2 + (\zeta + D - \frac{1}{2} D \beta) \cdot R^2 - \kappa r^1,\\ 
m_2 &= \zeta \cdot L^1 - \kappa l^2 \quad & s_2 = -\frac{1}{2} D\alpha \cdot L^1 + (\zeta + D - \frac{1}{2} D \alpha) \cdot R^1 - \kappa r^2.
\end{align*}

It is again a straight forward computation to check that
\begin{equation}\label{eq:zeta_mjsj}
(- \Delta - 2i \zeta \cdot \nabla + \omega^2 \varepsilon_0 \mu_0 + q_j)(m_j + s_j) = 0, \qquad j=1, 2.
\end{equation}
Further, using (\ref{ter:Lform_Rform}), (\ref{es:decay_R}), (\ref{es:decay_PR}) and
\[ \mathrm{supp}\, D\alpha \subset \overline{B(O;\rho)},\qquad \mathrm{supp}\, D\beta \subset \overline{B(O;\rho)}; \]
one sees that $ s_j \in L^2_\delta $. Recall that equation (\ref{eq:zeta_mjsj}) has a unique solution in $ L^2_\delta $ whenever $ |\zeta| > C(\delta, \rho) \norm{\omega^2 \varepsilon_0 \mu_0 + q_j}{}{L^\infty(B_\rho)} $, certainly it has to be the trivial one. Therefore, if $ m_j $ were also in $ L^2_\delta $, then $ f^j $ would vanish. Note that in order to have $ m_j \in L^2_\delta $, it is enough for its support to be compact. This could be accomplished by choosing $ l^j, L^j $ such that
\[ \zeta \cdot L^2 = \omega \varepsilon_0^{1/2} \mu_0^{1/2} l^1, \qquad \zeta \cdot L^1 = \omega \varepsilon_0^{1/2} \mu_0^{1/2} l^2. \]
\end{proof}

Next, we construct the same kind of solutions for the equation $ (P + W^*)\hat Y = 0 $. Since (\ref{eq:schr3}) holds and $ \omega^2 \varepsilon_0 \mu_0 + \hat Q $ has compact support in $ \overline{B(O; \rho)} $, the same kind of arguments used in the proof of Lemma \ref{le:CGO-sch} can be carried out to state the following lemma.
\begin{lemma}\label{le:CGO-rescale*} \rm
Let $ \delta $ be a constant such that $ -1 < \delta < 0 $ and let $ \zeta \in \mathbb{C}^3 $ be such that $ \zeta \cdot \zeta = \omega^2 \varepsilon_0 \mu_0 $ with
\[ |\zeta| > C \sum_{j,k=1}^8 \norm{(\omega^2 \varepsilon_0 \mu_0 I_8 + \hat Q)_j^k}{}{L^\infty(B_\rho)} \]
for some constant $ C = C(\delta, \rho) > 0 $.
Then, there exists a
\[ \hat Z = e^{i \zeta \cdot x} (\hat L + \hat R) \]
solution of $ (-\Delta I_8 + \hat Q) \hat Z = 0 $ in $ \mathbb{E} $, with $ \hat Z|_\inte \in H^2(\inte; \mathcal{Y}) $, $ \hat L = \hat L(\zeta) $ constant and $ \hat R = \hat R(\zeta) $ satisfying
\begin{align*}
\norm{\hat R}{}{L^2_\delta \mathcal{Y}} &\leq \frac{C(\delta, \rho)}{|\zeta|} |\hat L| \sum_{j,k=1}^8 \norm{(\omega^2 \varepsilon_0 \mu_0 I_8 + \hat Q)_j^k}{}{L^\infty(B_\rho)}\\
\norm{P \hat R}{}{L^2_\delta \mathcal{Y}} &\leq C(\delta, \rho) \sum_{j,k=1}^8 \norm{(\omega^2 \varepsilon_0 \mu_0 I_8 + \hat Q)_j^k}{}{L^\infty(B_\rho)} \Big(|\hat L| + \norm{\hat R}{}{L^2_\delta \mathcal{Y}} \Big).
\end{align*}
\end{lemma}

As a consequence of this lemma we get the following proposition, whose main difference with respect to Proposition \ref{pro:CGO-sch} is the construction of solutions for the rescaled system instead of solutions for Maxwell's equations.
\begin{proposition}\label{pro:CGO-rescale*}\rm
Let $ \zeta \in \mathbb{C}^3 $ be such that $ \zeta \cdot \zeta = \omega^2 \varepsilon_0 \mu_0 $ with
\[ |\zeta| > C(\rho) \sum_{j,k=1}^8 \norm{(\omega^2 \varepsilon_0 \mu_0 I_8 + \hat Q)_j^k}{}{L^\infty(B_\rho)} .\]
Then, there exists a
\[ \hat Y = e^{i\zeta \cdot x}(M + S) \]
solution for the equation $ (P + W^*) \hat Y = 0 $ in $ \mathbb{E} $, with $ \hat Y|_\inte \in H^1(\inte; \mathcal{Y}) $,
\[ M = \frac{1}{|\zeta|}\frv{\inner{\zeta}{\hat a}}{-\cross{\zeta}{\hat a}}{\inner{\zeta}{\hat b}}{\cross{\zeta}{\hat b}}, \]
for $ \hat a, \hat b $ constant complex vector fields, and $ S $ satisfying
\begin{equation*}
\norm{S}{}{L^2(\inte; \mathcal{Y})} \leq \frac{C(\rho, \inte)}{|\zeta|} \sum_{j,k=1}^8 \left( \norm{(\omega^2 \varepsilon_0 \mu_0 I_8 + \hat Q)_j^k}{}{L^\infty(B_\rho)} + \norm{W_j^k}{}{L^\infty(\inte)} \right).
\end{equation*}
\end{proposition}

\begin{proof}
Let $ \hat Z $ be the solution constructed in Lemma \ref{le:CGO-rescale*}, then by equation (\ref{eq:schr3}), $ \hat Y = (P - \overline{W}) \hat Z $ is a solution of $ (P + W^*) \hat Y = 0 $ in $ \mathbb{E} $. Considering
\[ \hat L = \frac{1}{|\zeta|} \fcv{0}{\hat b^t}{0}{\hat a^t}^t \]
with $ \hat a, \hat b $ constant complex vector fields, the solution $ \hat Y $ can be expressed as $ \hat Y = e^{i\zeta \cdot x}(M + S), $
with
\[ M = i P_\zeta \hat L, \qquad S = P \hat R + i P_\zeta \hat R - \overline{W} \hat L - \overline{W} \hat R, \]
where $ P_\zeta $ is as in (\ref{for:P_A}) and $ S $ satisfies
\begin{align*}
\norm{S}{}{L^2(\inte; \mathcal{Y})} \leq & C(\delta, \rho, \inte) \sum_{j,k=1}^8 \norm{(\omega^2 \varepsilon_0 \mu_0 I_8 + \hat Q)_j^k}{}{L^\infty(B_\rho)} \Big(|\hat L| + \norm{\hat R}{}{L^2_\delta \mathcal{Y}} \Big)\\
&+ C(\delta, \rho, \inte) |\hat L| \sum_{j,k=1}^8 \norm{(\omega^2 \varepsilon_0 \mu_0 I_8 + \hat Q)_j^k}{}{L^\infty(B_\rho)}\\
&+ \sum_{j,k=1}^8 \norm{W_j^k}{}{L^\infty(\inte)}\Big(\norm{\hat L}{}{L^2(\inte; \mathcal{Y})} + C(\inte) \norm{\hat R}{}{L^2_\delta \mathcal{Y}} \Big)\\
\leq & \frac{C(\delta, \rho, \inte)}{|\zeta|} \sum_{j,k=1}^8 \left( \norm{(\omega^2 \varepsilon_0 \mu_0 I_8 + \hat Q)_j^k}{}{L^\infty(B_\rho)} + \norm{W_j^k}{}{L^\infty(\inte)} \right).
\end{align*}
The last estimate is a consequence of Lemma \ref{le:CGO-rescale*}.
\end{proof}

\subsection{Proof of the log-type estimate}
The general ideas of this section go back to \cite{A}. The most relevant difference is the use of a Carleman estimate.

Let $ \mu_1, \gamma_1 $ and $ \mu_2, \gamma_2 $ be two pairs of coefficients under the hypothesis of Theorem \ref{th:stability} and let us choose
\begin{gather}
\zeta_1=-\frac{1}{2}\xi+i\left(\tau^2+\frac{|\xi|^2}{4}\right)^{1/2}\eta_1+ \left(\tau^2+ \omega^2 \varepsilon_0 \mu_0\right)^{1/2}\eta_2 \label{zeta1def}, \\
\zeta_2=\frac{1}{2}\xi-i\left(\tau^2+\frac{|\xi|^2}{4}\right)^{1/2}\eta_1+ \left(\tau^2+ \omega^2 \varepsilon_0 \mu_0\right)^{1/2}\eta_2 \label{zeta2def},
\end{gather}
with $ \tau \geq 1 $ a free parameter controlling the size of $ |\zeta_1| $ and $ |\zeta_2| $, where $ \xi, \eta_1,\eta_2 $ are constant vector fields satisfying $ |\eta_1| = |\eta_2| = 1 $, $ \eta_1 \cdot \eta_2 = 0 $ and $ \eta_j \cdot \xi = 0 $ for $ j = 1, 2 $. Note that $ \zeta_1 - \overline{\zeta_2} = - \xi $ and
\begin{equation*}
\frac{\zeta_1}{|\zeta_1|} = i \frac{\eta_1}{\sqrt{2}} + \frac{\eta_2}{\sqrt{2}} + \mathcal{O}(\tau^{-1}), \qquad \frac{\zeta_2}{|\zeta_2|} = -i \frac{\eta_1}{\sqrt{2}} + \frac{\eta_2}{\sqrt{2}} + \mathcal{O}(\tau^{-1}).
\end{equation*}
Let us consider $ Z_1 = e^{i \zeta_1 \cdot x} (L_1 + R_1), Y_1 $ the solutions stated in Proposition \ref{pro:CGO-sch} corresponding to the pair $ \mu_1, \gamma_1 $ with $ |\zeta_1| > C(\rho, M) $. Recall that
\begin{equation*}
L_1 = \frac{1}{|\zeta_1|} \frv{\zeta_1 \cdot a_1}{\omega \varepsilon_0^{1/2} \mu_0^{1/2} b_1}{\zeta_1 \cdot b_1}{\omega \varepsilon_0^{1/2} \mu_0^{1/2} a_1},\qquad
\norm{R_1}{}{L^2 (\inte; \mathcal{Y})} \leq \frac{C(\rho, \inte, M)}{|\zeta_1|}.
\end{equation*}
Additionally, consider $ Y_2 = e^{i \zeta_2 \cdot x} (M_2 + S_2) $ the solution stated in Proposition \ref{pro:CGO-rescale*} corresponding to $ \mu_2, \gamma_2 $ with $ |\zeta_2| > C(\rho, M) $. Also recall that
\begin{equation*}
M_2 = \frac{1}{|\zeta_2|} \frv{\zeta_2 \cdot a_2}{\cross{-\zeta_2}{a_2}}{\zeta_2 \cdot b_2}{\cross{\zeta_2}{b_2}},\qquad
\norm{S_2}{}{L^2(\inte; \mathcal{Y})} \leq \frac{C(\rho, \inte, M)}{|\zeta_2|}.
\end{equation*}
Next we plug these solutions into the estimate (\ref{es:suitablest}) of Proposition \ref{le:2-suitable-estimate}, with different choices of $ a_j, b_j $.

Choosing $ b_1 = b_2 = 0 $ and $ a_1, a_2 $ such that
\[ \left( i \frac{\eta_1}{\sqrt{2}} + \frac{\eta_2}{\sqrt{2}} \right) \cdot a_1 = \left( i \frac{\eta_1}{\sqrt{2}} + \frac{\eta_2}{\sqrt{2}} \right) \cdot \overline{a_2} = 1 \]
one gets, when $ \tau $ becomes large, that
\begin{gather*}
\Inner{(Q_1 - Q_2) Z_1}{Y_2}_\inte =\\
= \int_\inte e^{- i \xi \cdot x} \left( \frac{1}{2} \Delta (\alpha_1 - \alpha_2) + \frac{1}{4} ( \inner{\nabla \alpha_1}{\nabla \alpha_1} - \inner{\nabla \alpha_2}{\nabla \alpha_2} ) + ( \kappa^2_2 - \kappa^2_1 ) \right) \, dV\\
+\, \mathcal{O}((\tau^2 + |\xi|^2)^{-1/2}),
\end{gather*}
where the implicit constant in the symbol $ \mathcal{O} $ is $ C = C(\rho, \inte, M) $. Choosing $ a_1 = a_2 = 0 $ and $ b_1, b_2 $ such that
\[ \left( i \frac{\eta_1}{\sqrt{2}} + \frac{\eta_2}{\sqrt{2}} \right) \cdot b_1 = \left( i \frac{\eta_1}{\sqrt{2}} + \frac{\eta_2}{\sqrt{2}} \right) \cdot \overline{b_2} = 1 \]
one gets, when $ \tau $ becomes large, that
\begin{gather*}
\Inner{(Q_1 - Q_2) Z_1}{Y_2}_\inte =\\
= \int_\inte e^{- i \xi \cdot x} \left( \frac{1}{2} \Delta (\beta_1 - \beta_2) + \frac{1}{4} ( \inner{\nabla \beta_1}{\nabla \beta_1} - \inner{\nabla \beta_2}{\nabla \beta_2} ) + ( \kappa^2_2 - \kappa^2_1 ) \right) \, dV\\
+\, \mathcal{O}((\tau^2 + |\xi|^2)^{-1/2}),
\end{gather*}
where the implicit constant is $ C(\rho, \inte, M) $. Write
\begin{align*}
f &= \mathbf{1}_\inte \left( \frac{1}{2} \Delta (\alpha_1 - \alpha_2) + \frac{1}{4} ( \inner{\nabla \alpha_1}{\nabla \alpha_1} - \inner{\nabla \alpha_2}{\nabla \alpha_2} ) + ( \kappa^2_2 - \kappa^2_1 ) \right)\\
g &= \mathbf{1}_\inte \left( \frac{1}{2} \Delta (\beta_1 - \beta_2) + \frac{1}{4} ( \inner{\nabla \beta_1}{\nabla \beta_1} - \inner{\nabla \beta_2}{\nabla \beta_2} ) + ( \kappa^2_2 - \kappa^2_1 ) \right),
\end{align*}
where $ \mathbf{1}_\inte $ is the indicator function of $ \inte $. By Proposition \ref{le:2-suitable-estimate} and the properties of the special solutions, there exist three constants $ c = c(\inte) $, $ C = C(\rho, \inte, M) $ and $ C' = C'(\rho, M) $ such that, for any $ \tau \geq C' $ one has
\begin{equation*}
|\widehat{f}(\xi)| + |\widehat{g}(\xi)| \leq C \left( B\big( \delta_C(C_1, C_2) \big) e^{c(\tau^2 + |\xi|^2)^{1/2}} + (\tau^2 + |\xi|^2)^{-1/2} \right).
\end{equation*}

Note that, for $ s_1 < 0 $ and $ R \geq 1 $, one has
\begin{gather*}
\norm{f}{2}{H^{s_1}(\mathbb{E})} + \norm{g}{2}{H^{s_1}(\mathbb{E})} = \int_{|\xi| < R} (1 + |\xi|^2)^{s_1} \big( |\widehat{f}(\xi)|^2 + |\widehat{g}(\xi)|^2 \big) \, d\xi\\
+ \int_{|\xi| \geq R} (1 + |\xi|^2)^{s_1} \big( |\widehat{f}(\xi)|^2 + |\widehat{g}(\xi)|^2 \big) \, d\xi\\
\leq C \left( B\big( \delta_C(C_1, C_2) \big) e^{c(R + \tau)}+\, \tau^{-1} \right)^2 \int_0^R (1 + |r|^2)^{s_1} r^2 \, dr\\
+\, (1 + R^2)^{s_1} \left( \norm{f}{2}{L^2(\inte)} + \norm{g}{2}{L^2(\inte)} \right).
\end{gather*}
Therefore,
\begin{equation*}
\norm{f}{}{H^{s_1}(\mathbb{E})} + \norm{g}{}{H^{s_1}(\mathbb{E})} \leq C \left( B\big( \delta_C(C_1, C_2) \big) e^{c(R + \tau)} + \tau^{-1}  R^{3/2 + s_1} + R^{s_1} \right).
\end{equation*}
Now we choose $ R $ in such a way that $ \tau^{-1}  R^{3/2 + s_1} $ decays as $ R^{s_1} $, that is, $ R = \tau^{2/3} $, hence
\begin{equation*}
\norm{f}{}{H^{s_1}_0(\inte)} + \norm{g}{}{H^{s_1}_0(\inte)} \leq C \left( B\big( \delta_C(C_1, C_2) \big) e^{c \tau} + \tau^{2/3 s_1} \right).
\end{equation*}
On the other hand, the a priori bound was chosen to have
\[ \norm{f}{}{H^{s_2}(\inte)} + \norm{g}{}{H^{s_2}(\inte)} \leq C(M), \]
for $ 0 < s_2 < 1/2 $. Finally, by the interpolation estimate (\ref{es:interpolation}) there exist two constants $ C' = C'(\rho, M) $ and $ C = C(\rho, \inte, M, \omega) $ such that, for any $ \tau \geq C' $, the following estimate holds
\begin{equation}\label{es:normL2de_f_g}
\norm{f}{}{L^2(\inte)} + \norm{g}{}{L^2(\inte)} \leq C \left( B\big( \delta_C(C_1, C_2) \big) e^{c \tau} + \tau^{2/3 s_1} \right)^\theta,
\end{equation}
with $ 0 = \theta s_1 + (1 - \theta) s_2 $.

The idea now is to transfer this estimate from $ f, g $ to the difference of the coefficients $ \mu_1 - \mu_2 $ and $ \gamma_1 - \gamma_2 $. This can be accomplished by using the following Carleman estimate.

There exists a positive constant $ C(\inte) $ such that, for all $ h \leq 1 $ and any function $ \phi \in C^{1,1}(\overline{\inte}) $, the following estimate holds
\begin{gather*}
h \norm{e^{\varphi / h}\phi}{2}{L^2(\inte)} + h^3 \norm{e^{\varphi / h} \nabla \phi}{2}{L^2(\inte; \mathbb{C}^3)} \leq\\
\leq C \left( h^4 \norm{e^{\varphi / h} \Delta \phi}{2}{L^2(\inte)} + h \norm{e^{\varphi / h}\phi}{2}{L^2(\bou)} + h^3 \norm{e^{\varphi / h} \nabla \phi}{2}{L^2(\bou; \mathbb{C}^3)} \right),
\end{gather*}
where $ \varphi = 1/2 |x-x_0|^2 $ with $ x_0 \notin \overline{\inte} $. The constant here depends on the distance from $ x_0 $ to $ \inte $ and on the diameter of $ \inte $. A Carleman estimate of this type can be found in \cite{I}.

A simple computation give:
\begin{align*}
f &= \mathbf{1}_\inte \gamma_1^{-1/2} \left[ \Delta(\gamma_1^{1/2} - \gamma_2^{1/2}) + q_f(\gamma_1^{1/2} - \gamma_2^{1/2}) + p_f(\mu_1^{1/2} - \mu_2^{1/2}) \right],\\
g &= \mathbf{1}_\inte \mu_1^{-1/2} \left[ \Delta(\mu_1^{1/2} - \mu_2^{1/2}) + q_g(\mu_1^{1/2} - \mu_2^{1/2}) + p_g(\gamma_1^{1/2} - \gamma_2^{1/2}) \right];
\end{align*}
where
\begin{align*}
q_f = - \left( \frac{\Delta \gamma_2^{1/2}}{\gamma_2^{1/2}} + \omega^2 \gamma_1^{1/2}(\gamma_1^{1/2} \mu_1 + \gamma_2^{1/2} \mu_2) \right),& &
p_f = - \omega^2 \gamma_1 \gamma_2^{1/2} (\mu_1^{1/2} + \mu_2^{1/2}),\\
q_g = - \left( \frac{\Delta \mu_2^{1/2}}{\mu_2^{1/2}} + \omega^2 \mu_1^{1/2}(\mu_1^{1/2} \gamma_1 + \mu_2^{1/2} \gamma_2) \right),& &
p_g = - \omega^2 \mu_1 \mu_2^{1/2} (\gamma_1^{1/2} + \gamma_2^{1/2}).
\end{align*}
Note that, thanks to the a priori bounds, we have the following differential inequalities:
\begin{align*}
|\Delta(\gamma_1^{1/2} - \gamma_2^{1/2})| \leq C(M) (|f| + |\gamma_1^{1/2} - \gamma_2^{1/2}| + |\mu_1^{1/2} - \mu_2^{1/2}|),\\
|\Delta(\mu_1^{1/2} - \mu_2^{1/2})| \leq C(M) (|g| + |\gamma_1^{1/2} - \gamma_2^{1/2}| + |\mu_1^{1/2} - \mu_2^{1/2}|).
\end{align*}

In order to simplify the notation, let us write $ \phi_1 = \gamma_1^{1/2} - \gamma_2^{1/2} $ and $ \phi_2 = \mu_1^{1/2} - \mu_2^{1/2} $. By the differential inequalities written above and the Carleman estimate, one has
\begin{gather*}
\sum_{j=1,2} \left( h \norm{e^{\varphi / h}\phi_j}{2}{L^2(\inte)} + h^3 \norm{e^{\varphi / h} \nabla \phi_j}{2}{L^2(\inte; \mathbb{C}^3)} \right) \leq\\
\leq C'' \sum_{j=1,2} \left( h^4 \norm{e^{\varphi / h} \phi_j}{2}{L^2(\inte)} + h \norm{e^{\varphi / h} \phi_j }{2}{L^2(\bou)} + h^3 \norm{e^{\varphi / h} \nabla \phi_j}{2}{L^2(\bou; \mathbb{C}^3)} \right)\\
+ C'' h^4 \left( \norm{e^{\varphi / h} f}{2}{L^2(\inte)} + \norm{e^{\varphi / h} g}{2}{L^2(\inte)} \right),
\end{gather*}
where the constant is $ C'' = C''(\inte, M) $. The terms $ h^4 \norm{e^{\varphi / h} \phi_j}{2}{L^2(\inte)} $, with $ j = 1,2 $, can be absorbed by the left hand side of the inequality. Hence, if $ d_1 = \textrm{inf} \{ d_e(x; x_0)^2 : x \in \inte \} $ and $ d_2 = \textrm{sup} \{ d_e(x; x_0)^2 : x \in \inte \} $ we get, for any $ h < C''(\inte, M)^{-1/3} $, that
\begin{gather*}
e^{d_1 / h} \sum_{j=1,2} \left( h \norm{\phi_j}{2}{L^2(\inte)} + h^3 \norm{\nabla \phi_j}{2}{L^2(\inte; \mathbb{C}^3)} \right) \leq C'' e^{d_2 / h}\times\\
\left[ h^4 \left( \norm{f}{2}{L^2(\inte)} + \norm{g}{2}{L^2(\inte)} \right) + \sum_{j=1,2} \left( h \norm{\phi_j }{2}{L^2(\bou)} + h^3 \norm{\nabla \phi_j}{2}{L^2(\bou; \mathbb{C}^3)} \right)\right ].
\end{gather*}
But now we can easily estimate
\begin{gather*}
\norm{\phi_1}{}{L^2(\bou)} \leq C \norm{\gamma_1 - \gamma_2}{}{L^\infty(\bou)} \leq C B \big( \delta_C(C_1, C_2) \big),\\
\begin{align*}
\norm{\nabla \phi_1}{}{L^2(\bou; \mathbb{C}^3)} &\leq C \left( \norm{\gamma_1 - \gamma_2}{}{L^\infty(\bou)} + \norm{\nabla (\gamma_1 - \gamma_2)}{}{L^\infty(\bou; \mathbb{C}^3)} \right)\\
&\leq C B \big( \delta_C(C_1, C_2) \big),
\end{align*}\\
\norm{\phi_2}{}{L^2(\bou)} \leq C \norm{\mu_1 - \mu_2}{}{L^\infty(\bou)} \leq C B \big( \delta_C(C_1, C_2) \big),\\
\begin{align*}
\norm{\nabla \phi_2}{}{L^2(\bou; \mathbb{C}^3)} &\leq C \left( \norm{\mu_1 - \mu_2}{}{L^\infty(\bou)} + \norm{\nabla(\mu_1 - \mu_2)}{}{L^\infty(\bou; \mathbb{C}^3)} \right)\\
&\leq C B \big( \delta_C(C_1, C_2) \big),
\end{align*}\\
\norm{\gamma_1 - \gamma_2}{}{L^2(\inte)} + \norm{\nabla (\gamma_1 - \gamma_2)}{}{L^2(\inte; \mathbb{C}^3)} \leq C \left( \norm{\phi_1}{}{L^2(\inte)} + \norm{\nabla \phi_1}{}{L^2(\inte; \mathbb{C}^3)} \right),\\
\norm{\mu_1 - \mu_2}{}{L^2(\inte)} + \norm{\nabla (\mu_1 - \mu_2)}{}{L^2(\inte; \mathbb{C}^3)} \leq C \left( \norm{\phi_2}{}{L^2(\inte)} + \norm{\nabla \phi_2}{}{L^2(\inte; \mathbb{C}^3)} \right).
\end{gather*}
The constants above depend on the a priori bounds $ M $. These inequalities and estimate (\ref{es:normL2de_f_g}) gives us
\begin{align*}
\norm{\gamma_1 - \gamma_2}{}{H^1(\inte)} + \norm{\mu_1 - \mu_2}{}{H^1(\inte)} \leq & C e^{\frac{d_2 - d_1}{2h}} \left( B\big( \delta_C(C_1, C_2) \big) e^{c \tau} + \tau^{2/3 s_1} \right)^\frac{s_2}{s_2 - s_1}\\
& + C e^{\frac{d_2 - d_1}{2h}} B\big( \delta_C(C_1, C_2) \big),
\end{align*}
where $ d_2 > d_1 $, $ s_1 < 0 < s_2 < 1/2 $, $ c = c(\inte) $, $ C = C(\rho, \inte, M) $, $ \tau \geq C'(\rho, M) $ and $ h < C''(\inte, M)^{-1/3} $. To end up with the estimate given in the statement, it is enough to note that
\[ 0 < -\frac{2}{3}\frac{s_1 s_2}{s_2 - s_1} < \frac{2}{3} s_2, \]
and to choose the parameter $ \tau $ as
\[ \tau = - \frac{1}{2c} \mathrm{log} \, B \big( \delta_C(C_1, C_2) \big). \]


\appendix
\section{Appendix}
In this appendix, we state some elementary concepts and collect some known facts relative to Lipschitz domains and Sobolev-Besov spaces.

\subsection{Lipschitz domain}
\begin{definition}\label{def:Lipschitz-domain}\rm
Let $ \inte $ be a nonempty, proper open subset of $ \mathbb{E} $ and fix a point $ x_0 $ on its boundary $ \bou $. We say that $ \inte $ is a \emph{Lipschitz domain near} $ x_0 $ if there exist 
\begin{itemize}
\item[(i)] a plane $ q \subset \mathbb{E} $ passing through $ x_0 $ and a choice of a unit vector $ N_q $ normal  to $ q $;
\item[(ii)] some euclidean coordinates $ \mathcal{E}_0: \mathbb{E} \rightarrow \mathbb{R}^3 $ such that $ \mathcal{E}_0(x_0)^j=0 $ for $ j= 1,2,3 $ and $ \mathcal{E}_0(x) \in \mathbb{R}^2 \times\{0\} $, for any $ x \in q $ (for short, we shall denote $ \mathcal{E}_0(x)^j $ by $ y^j $);
\item[(iii)] and an open cylinder $ C^{x_0}_{c_1,c_2} = \{x \in \mathbb{E}: |y'|<c_1, |y^3|<c_2\} $ --called coordinate cylinder near $ x_0 $-- such that
\begin{gather*}
C^{x_0}_{c_1,c_2} \cap \inte = C^{x_0}_{c_1,c_2} \cap \{x \in \mathbb{E}: y^3 > \phi(y^1, y^2)\},\\
C^{x_0}_{c_1,c_2} \cap \bou = C^{x_0}_{c_1,c_2} \cap \{x \in \mathbb{E}: y^3 = \phi(y^1, y^2)\},\\
C^{x_0}_{c_1,c_2} \cap \overline{\inte}^c = C^{x_0}_{c_1,c_2} \cap \{x\in\mathbb{E}: y^3 < \phi(y^1, y^2)\};
\end{gather*}
for some Lipschitz function $ \phi: \mathbb{R}^2\rightarrow \mathbb{R} $ satisfying
\[ \phi(0) = 0,\quad \textrm{and}\quad |\phi(y^1, y^2)| < c_2 \quad \textrm{if}\quad |y'| \leq c_1. \]
Finally, we say that $ \inte $ is a \emph{Lipschitz domain} if it is a Lipschitz domain near every point $ x \in \bou $.
\end{itemize}
\end{definition}
In this definition the superscript $ c $ denotes the complement of a set, relative to $ \mathbb{E} $; and $ |y'|^2 = |y^1|^2 + |y^2|^2 $.

Recall that for a Lipschitz domain there exists a measurable unit normal vector field $ N $ along $ \bou $ pointing outward. In this context, we can set the following integration by parts formulas
\begin{equation}\label{eq:delta-d_identity}
- \int_\inte (\diver u) f \,dV = \int_\inte \inner{u}{\nabla f} \,dV - \int_{\bou} (\inner{N}{u}) f|_{\bou} \,dA
\end{equation}
and
\begin{equation}\label{eq:*d_identity}
\int_\inte \inner{(\curl u)}{v} \,dV = \int_\inte \inner{u}{(\curl v)} \,dV + \int_{\bou} \inner{(\cross{N}{u})}{v|_{\bou}} \,dA,
\end{equation}
where $ f $ is a smooth function on $ \overline{\inte} $ and $ u, v $ are vector fields on $ \overline{\inte} $. Here $ dA $ stands for the \emph{area element} and $ \centerdot|_{\bou} $ denotes the restriction to $ \bou $.

\subsection{Sobolev and Besov spaces}
Most of the facts collected here can be found in \cite{JK} and the references therein.
\begin{definition}\rm
For any $ s\in\mathbb{R} $, define the potential Sobolev space on $ \mathbb{E} $ as
\[ H^s(\mathbb{E}) = \overline{C^\infty_0(\mathbb{E})}^{\norm{.}{}{H^s(\mathbb{E})}}, \]
with the norm
\[ \norm{f}{2}{H^s(\mathbb{E})}  = \int_{\mathbb{R}^3} (1+|\xi|^2)^s |\widehat{f}(\xi)|^2 \, d\xi.\]
When $ s = 0 $ this space is $ L^ 2(\mathbb{E}) $.
\end{definition}

Denote by $f|_\inte$ the restriction of $f\in L^2(\mathbb{E})$ to $\inte\subset\mathbb{E}$.
\begin{definition}\rm
For any $ s \geq 0 $, define the potential Sobolev space on $ \inte $ as
\[ H^s(\inte) = \{ f|_\inte : f \in H^s(\mathbb{E}) \}, \]
with the norm
\[ \norm{g}{}{H^s(\inte)} = \textrm{inf} \{ \norm{f}{}{H^s(\mathbb{E})} : f|_\inte = g \}. \]
On the other hand, for $ s \in \mathbb{R} $, define
\[ H^s_0(\inte) = \{ f \in H^s(\mathbb{E}) : \textrm{supp}\,f \subseteq \overline{\inte} \}, \]
with norm
\[ \norm{f}{}{H^s_0(\inte)} = \norm{f}{}{H^s(\mathbb{E})}. \]
Finally, for $ s > 0 $, define the space $ H^{-s}(\inte) $ as the dual of $ H^s_0(\inte) $, that is,
\[ H^{-s}(\inte) = (H^s_0(\inte))^*. \]
It is well known that $ C^{\infty}(\overline{\inte}) = \{ f|_{\overline{\inte}} : f \in C^{\infty}(\mathbb{E}) \} $ is dense in $ H^s(\inte) $, for $ s \in \mathbb{R} $. When $ s = 0 $ this is $ L^2(\inte) $.
\end{definition}

On the other hand, $ H^s(\mathbb{E}) $ is a complex interpolation scale for $ s \in \mathbb{R} $; that is, for any $ s_1, s_2 \in \mathbb{R} $, one has
\[ [ H^{s_1}(\mathbb{E}) , H^{s_2}(\mathbb{E}) ]_\theta = H^s(\mathbb{E}), \]
with $ s = \theta s_1 + (1-\theta) s_2 $ with $ \theta \in (0,1)$. Moreover, the extension by zero outside $ \inte $ is a bounded linear operator from $ H^s(\inte)$ to $ H^s(\mathbb{E}) $ for $ -1/2 < s < 1/2$ (see \cite{T}), which allows to identify the spaces $ H^s(\inte) $ and $ H^s_0(\inte) $ for the same range of $ s $. With these facts in mind, notice that
\begin{equation}\label{es:interpolation}
\norm{f}{}{H^s_0(\inte)} \leq C \norm{f}{\theta}{H^{s_1}_0(\inte)} \norm{f}{1-\theta}{H^{s_2}(\inte)},
\end{equation}
for $ s = \theta s_1 + (1-\theta) s_2 $ with $ s_1 \in \mathbb{R} $ and $ -1/2 < s_2 < 1/2 $.

Recall that
\[ W^{2,\infty}(\inte) = \{ f \in L^\infty(\inte) : \partial^\alpha f \in L^\infty(\inte),\, 0 < |\alpha|\leq 2 \}, \]
with the norm
\[ \norm{f}{}{W^{2,\infty}(\inte)} = \sum_{0 \leq |\alpha| \leq 2} \norm{\partial^\alpha f}{}{L^\infty(\inte)}. \]

We next define the Besov spaces in an intrinsic way. To do it we introduce the functional
\[ I_s(f) = \left( \int_{\mathbb{R}^2} \frac{\norm{f(\centerdot + y)-f(\centerdot)}{2}{L^2(\mathbb{R}^2)}}{|y|^{2+2s}} \,dy \right)^{1/2}, \]
defined for $ f \in \mathcal{S}(\mathbb{R}^2) $ --the space of rapidly decreasing functions.
\begin{definition}\rm
For $ 0 < s < 1 $, let us define the Besov space
\[ B^s(\mathbb{R}^2) = \{ f \in L^2(\mathbb{R}^2) : I^s(f) < + \infty \}, \]
with the norm
\[ \norm{f}{}{B^s(\mathbb{R}^2)} = \norm{f}{}{L^2(\mathbb{R}^2)} + I^s(f). \]
\end{definition}

Now we shall extend these Besov spaces on $ \mathbb{R}^2 $ to Besov spaces on $ \bou $. Let $ x_1, \dots, x_n $ belong to $ \bou $ and $ \Gamma_1, \dots, \Gamma_n $ be $ \Gamma_j = C^{x_j}_{c_1,c_2} \cap \bou $ for $ j=1, \dots, n $, such that $ \bou = \Gamma_1 \cup \dots \cup \Gamma_n $; and consider a partition of unity $ \chi_1, \dots, \chi_n $ subordinate to $ \Gamma_1, \dots, \Gamma_n $. We shall say that $f \in B^s(\bou) $ for $ 0 < s < 1 $ if
\[ (\chi_j f) \circ \mathcal{E}^{-1}_j(\centerdot, \phi_j(\centerdot)) \in B^s(\mathbb{R}^2), \]
for any possible choice of points and any partition of unity related to them as above. Here $ \mathcal{E}_j $ and  $ \phi_j $ are, respectively, the euclidean coordinates and the function defining the boundary locally, corresponding to the point $ x_j $. The norm defined on these spaces will be given by
\begin{gather*}
\norm{f}{}{B^s(\bou)} = \textrm{inf}\, \{ \sum_{j= 1}^n \norm{(\chi_j f)\circ \mathcal{E}^{-1}_j(\centerdot,\phi_j(\centerdot))}{}{B^s(\mathbb{R}^2)} : n \in \mathbb{N},\\
x_j \in \bou,\, \textrm{supp}(\chi_j) \subset \Gamma_j, \, j=1, \dots ,n \}.
\end{gather*}

One of the reasons to introduce these spaces is to describe the properties of the trace operator. Namely, the trace operator
$ \centerdot|_{\bou} : H^s(\inte) \longrightarrow B^{s-1/2}(\bou) $ is well-defined, bounded and onto, whenever $ 1/2 < s < 3/2 $. In addition, it has a bounded right inverse whose norm is controlled by $ s $ and the Lipschitz character of $ \inte $.

For later references, we give the following lemma.
\begin{lemma}\label{le:besov_product}\rm
Let $ s, \epsilon $ be such that $ 0 < s < 1 $ and $ 0 < \epsilon \leq 1-s $. Then, there exists a constant $ C(s,\epsilon) > 0 $ such that, for any $ g \in C^{0,s + \epsilon}(\bou) $ and any $ f \in B^s(\bou) $,
\begin{equation}\label{es:besov_product}
\norm{gf}{}{B^s(\bou)} \leq C \norm{g}{}{C^{0,s+\epsilon}(\bou)} \norm{f}{}{B^s(\bou)}.
\end{equation}
\end{lemma}

\textbf{Remark:} The constant $ C(s,\epsilon) $ given here blows up when $ \epsilon $ becomes small.

\begin{proof} Let $ x_1, \dots, x_n $ belong to $ \bou $ and $ \Gamma_1, \dots, \Gamma_n $ be such that $ \Gamma_j = C^{x_j}_{c_1,c_2} \cap \bou $ for $ j = 1, \dots, n $, such that $ \bou = \Gamma_1 \cup \dots \cup \Gamma_n $. Let $ \lambda > 0 $ be the Lebesgue number associated to $ \{ \Gamma_j \}_{j = 1}^n $ and define
\[ \tilde{\Gamma}_j = \{ x \in \Gamma_j : \inf_{x' \in \bou \setminus \Gamma_j} d_e (x, x') > \lambda / 4 \}. \]
Note that $ \bou = \tilde{\Gamma}_1 \cup \dots \cup \tilde{\Gamma}_n $. Let us consider a partition of unity $ \chi_1, \dots, \chi_n $ subordinated to $ \tilde{\Gamma}_1, \dots, \tilde{\Gamma}_n $ and denote
\begin{align*}
f_j(y) = ( \chi_j f ) \circ \mathcal{E}^{-1}_j(y,\phi_j(y)),& &
g_j(y) = ( \mathbf{1}_{\Gamma_j} g ) \circ \mathcal{E}^{-1}_j(y,\phi_j(y)).
\end{align*}
Here $ \mathbf{1}_{\Gamma_j} $ stands for the indicator function of $ \Gamma_j $. Consider $ x_k \in \Gamma_j $ with $ k = 1,2 $, we can write $ x_k = \mathcal{E}^{-1}_j(y_k,\phi_j(y_k)) $ for $ y_k \in \mathbb{R}^2 $. By the Lipschitz character of $ \phi_j $, one has
\[ d_e(x_1, x_2) \leq C |y_1 - y_2|. \]
Then, noting that
\begin{gather*}
\left( I^s(g_j f_j) \right)^2 \leq  \, 2 \sup_{\substack{y' \in \mathcal{E}_j(\tilde{\Gamma}_j) \\ C|y| < \lambda / 5}} \frac{|g_j(y' + y)-g_j(y')|^2}{|y|^{2(s+\epsilon)}} \norm{f_j}{2}{L^2(\mathbb{R}^2)} \int_{C|y| < \lambda / 5} \frac{1}{|y|^{2(1-\epsilon)}} \, dy\\
+ 8 \norm{g_j}{2}{L^\infty(\mathbb{R}^2)} \norm{f_j}{2}{L^2(\mathbb{R}^2)} \int_{C|y| \geq \lambda / 5} \frac{1}{|y|^{2(1+s)}} \, dy + 2 \norm{g_j}{2}{L^\infty(\mathbb{R}^2)} \left( I^s(f_j) \right)^2,
\end{gather*}
we can achieve the result.
\end{proof}

\begin{definition}\rm
For $ 0 < s < 1 $, define the space $ B^{-s}(\bou) $ as the dual of $ B^s(\bou) $, that is,
\[ B^{-s}(\bou) = (B^s(\bou))^\ast. \]
\end{definition}
In the same conditions as in Lemma \ref{le:besov_product}, one has, by duality --recall that $ \dual{gf}{h} = \dual{f}{\overline{g}h} $ for any $ h \in B^s(\bou) $--, that
\[ \norm{gf}{}{B^{-s}(\bou)} \leq C \norm{g}{}{C^{0,s+\epsilon}(\bou)} \norm{f}{}{B^{-s}(\bou)}. \]
The constant here is the same as in Lemma \ref{le:besov_product} and blows up when $ \epsilon $ becomes small.

\begin{definition}\rm
For $ s \in \mathbb{R} $ define
\[ H^s(\inte; \mathbb{C}^3) = \overline{ \mathcal{X}\mathbb{E}|_{\overline{\inte}} }^{\norm{\centerdot}{}{H^s(\inte; \mathbb{C}^3)}} \]
where
\[ \norm{u}{2}{H^s(\inte;  \mathbb{C}^3)} = \sum_{j = 1}^3 \norm{u^{(j)}}{2}{H^s(\inte)}. \]
Here $ \mathcal{X}\mathbb{E}|_{\overline{\inte}} = \{ u|_{\overline{\inte}} : u \in \mathcal{X}\mathbb{E} \} $. When $ s = 0 $ this space will be denoted by $ L^2(\inte; \mathbb{C}^3) $.
\end{definition}

\begin{definition}\rm
For $ 0 < |s| < 1 $ define
\[ B^s(\bou; \mathbb{C}^3) = \overline{ \mathcal{X}\mathbb{E}|_{\bou} }^{\norm{\centerdot}{}{B^s(\bou; \mathbb{C}^3)}} \]
where
\[ \norm{w}{2}{B^s(\bou; \mathbb{C}^3)} = \sum_{j = 1}^3 \norm{w^{(j)}}{2}{B^s(\bou)} \]
and $ \mathcal{X}\mathbb{E}|_{\bou} = \{ u|_{\bou} : u \in \mathcal{X}\mathbb{E} \} $. In  the same way, let us define
\[ L^2(\bou; \mathbb{C}^3) = \overline{ \mathcal{X}\mathbb{E}|_{\bou} }^{\norm{\centerdot}{}{L^2(\bou; \mathbb{C}^3)}} \]
where
\[ \norm{w}{2}{L^2(\bou; \mathbb{C}^3)} = \sum_{j = 1}^3 \norm{w^{(j)}}{2}{L^2(\bou)}. \]
\end{definition}

Let us point out that, for $ 0< |s| < 1 $ and $ 0< \epsilon \leq 1-|s| $,
\[ \norm{g w}{}{B^s(\bou; \mathbb{C}^3)} \leq C \norm{g}{}{C^{0,|s|+\epsilon}(\bou)} \norm{w}{}{B^s(\bou; \mathbb{C}^3)}. \]
The constant above is the same as the one in Lemma \ref{le:besov_product} and, once again, it blows up when $ \epsilon $ becomes small.

On the other hand, by using the trace operator component by component
\[ \centerdot|_{\bou} : H^s(\inte; \mathbb{C}^3) \longrightarrow B^{s-1/2}(\bou; \mathbb{C}^3) \]
is bounded and onto, whenever $ 1/2 < s < 3/2 $.


\end{document}